\documentclass[journal]{IEEEtran}

\usepackage{cite}

\ifCLASSINFOpdf
\usepackage[pdftex]{graphicx}
\usepackage{epstopdf}
\DeclareGraphicsExtensions{.pdf,.jpeg,.png}
\else
\usepackage[dvips]{graphicx}
\DeclareGraphicsExtensions{.eps}
\fi

\usepackage[cmex10]{amsmath}

\usepackage{algorithmic}

\usepackage{array}
\usepackage{url}

\ifCLASSOPTIONcompsoc
\usepackage[caption=false,font=normalsize,labelfont=sf,textfont=sf]{subfig}
\else
\usepackage[caption=false,font=footnotesize]{subfig}
\fi

\usepackage{graphicx}

\usepackage[T1]{fontenc} % textbf
\usepackage{footnote}
\usepackage{multirow}
\usepackage{amsthm}

\PassOptionsToPackage{normalem}{ulem} % subrayado
\usepackage{ulem} % subrayado

\theoremstyle{plain}

\theoremstyle{plain}

\theoremstyle{definition}

\providecommand{\definitionname}{Definition}
\providecommand{\lemmaname}{Lemma}
\providecommand{\theoremname}{Theorem}

\usepackage{pgfplots}
\usepackage{pgfplotstable}
\pgfplotsset{compat=newest}

\usepackage{multirow}
\usepackage{rotating}

\begin{document}

\title{An Assessment on the Use of Stationary Vehicles as a Support to Cooperative Positioning\thanks{Some preliminary material was presented at the 3rd International Conference on Connected Vehicles and Expo, Wien, Austria, 2014.}}

\author{Rodrigo~H.~Ord\'{o}\~{n}ez-Hurtado,
Emanuele Crisostomi,
Wynita~M.~Griggs
and~Robert~N.~Shorten % <-this % stops a space

\thanks{R. Ord\'{o}\~{n}ez-Hurtado, W. Griggs and R. Shorten are with The Hamilton Institute, National University of Ireland Maynooth, Maynooth, Co. Kildare, Ireland.}% <-this % stops a space
\thanks{E. Crisostomi is with the Department of Energy, Systems, Territory and Construction Engineering, University of Pisa, Italy.}
\thanks{R. Shorten is also with IBM Research Ireland, Dublin, Ireland.}
}

\maketitle

\begin{abstract}
In this paper, we consider the use of stationary vehicles as tools to enhance the localisation capabilities of moving vehicles in a VANET. We examine the idea in terms of its potential benefits, technical requirements, algorithmic design and experimental evaluation. Simulation results are given to illustrate the efficacy of the technique.
\end{abstract}

\begin{IEEEkeywords}
Cooperative Positioning, V2V Communications, Stationary vehicles, Parked vehicles, queue.
\end{IEEEkeywords}

\IEEEpeerreviewmaketitle

%¶¶¶¶¶¶¶¶¶¶¶¶¶¶¶¶¶¶¶¶¶¶¶¶¶¶¶¶¶¶¶¶¶¶¶¶¶¶¶¶¶¶¶¶¶¶¶¶¶¶¶¶¶¶¶¶¶¶¶¶¶¶¶¶¶¶¶¶¶¶¶¶¶¶¶¶¶¶¶¶¶¶¶¶¶¶¶¶¶¶¶¶¶¶¶¶¶¶¶¶¶¶¶¶¶¶¶¶¶¶¶¶¶¶¶¶¶¶¶¶¶¶¶¶¶¶¶¶¶¶¶¶¶¶¶¶¶¶¶¶¶¶¶¶¶¶¶¶¶¶¶¶
%¶¶¶¶¶¶¶¶¶¶¶¶¶¶¶¶¶¶¶¶¶¶¶¶¶¶¶¶¶¶¶¶¶¶¶¶¶¶¶¶¶¶¶¶¶¶¶¶¶¶¶¶¶¶¶¶¶¶¶¶¶¶¶¶¶¶¶¶¶¶¶¶¶¶¶¶¶¶¶¶¶¶¶¶¶¶¶¶¶¶¶¶¶¶¶¶¶¶¶¶¶¶¶¶¶¶¶¶¶¶¶¶¶¶¶¶¶¶¶¶¶¶¶¶¶¶¶¶¶¶¶¶¶¶¶¶¶¶¶¶¶¶¶¶¶¶¶¶¶¶¶¶
\section{Introduction}
%¶¶¶¶¶¶¶¶¶¶¶¶¶¶¶¶¶¶¶¶¶¶¶¶¶¶¶¶¶¶¶¶¶¶¶¶¶¶¶¶¶¶¶¶¶¶¶¶¶¶¶¶¶¶¶¶¶¶¶¶¶¶¶¶¶¶¶¶¶¶¶¶¶¶¶¶¶¶¶¶¶¶¶¶¶¶¶¶¶¶¶¶¶¶¶¶¶¶¶¶¶¶¶¶¶¶¶¶¶¶¶¶¶¶¶¶¶¶¶¶¶¶¶¶¶¶¶¶¶¶¶¶¶¶¶¶¶¶¶¶¶¶¶¶¶¶¶¶¶¶¶¶
%¶¶¶¶¶¶¶¶¶¶¶¶¶¶¶¶¶¶¶¶¶¶¶¶¶¶¶¶¶¶¶¶¶¶¶¶¶¶¶¶¶¶¶¶¶¶¶¶¶¶¶¶¶¶¶¶¶¶¶¶¶¶¶¶¶¶¶¶¶¶¶¶¶¶¶¶¶¶¶¶¶¶¶¶¶¶¶¶¶¶¶¶¶¶¶¶¶¶¶¶¶¶¶¶¶¶¶¶¶¶¶¶¶¶¶¶¶¶¶¶¶¶¶¶¶¶¶¶¶¶¶¶¶¶¶¶¶¶¶¶¶¶¶¶¶¶¶¶¶¶¶¶

Cooperative positioning (CP) for vehicular networks is a very topical problem. Exact road positioning is viewed as a key enabler of services such as road pricing (lane pricing), lane prioritisation for special vehicles (electric vehicles), and is also recognised as a key technology to enable autonomous driving. In this context, CP is being viewed as a means of overcoming the shortcomings of traditional Global Navigation Satellite systems (GNSS) for vehicular applications.

There already exist a range of CP techniques for realising on-road localisation. Many of these techniques rely heavily on support from dedicated infrastructure.
Examples include Differential Global Positioning Systems (DGPS), Real Time Kinematic (RTK) positioning, and Assisted Global Positioning Systems (A-GPS).
There are many excellent papers on CP, and ideas for realising this technology have appeared in various areas. The interested reader is referred to any one of the excellent papers on this topic; in particular, see the recent survey paper \cite{Alam2013}.

The main contribution of this paper is to assess the potential of stationary vehicles to support the current CP systems in a Vehicular Ad-hoc Network (VANET).
This idea follows recent works that have recognised the ability of suitably equipped stationary vehicles, in the absence of any time constraints, to pinpoint their own locations in a precise manner, and have recommended their usage as general service delivery platforms \cite{ibmICCVE}. We shall show that such a novel idea has the potential of greatly improving current CP systems, without requiring any dedicated infrastructure, and at the same time opening a new market for vehicle owners and car manufacturers to monetise vehicles in a non-traditional manner. Our vision imagines that vehicle owners make their vehicles available for such applications in exchange for monetary compensation, or for access to privileged services and facilities (e.g., parking spots).

Our paper is organised as follows. In Section \ref{stateofart} of the paper, we examine the state-of-the art concerning the use of stationary vehicles in providing VANET-related services. In Section \ref{S:Characteristics}, we provide a discussion to assess the potential of stationary powered-on and powered-off vehicles in Intelligent Transportation Systems (ITS) applications. Section \ref{S:CP_approach} formally describes how stationary vehicles can be integrated in CP applications in practice, and some experimental evaluations are presented in Section \ref{S:Evaluation}. Finally, Section \ref{S:Conclusions} concludes the paper.

\section{State of the Art}\label{stateofart}

Cooperative positioning is a topic that is of interest in many communities. As we have mentioned, the interested reader is referred to the excellent survey covering CP technologies in ITS applications \cite{Alam2013}, and to the EU Funded Project TEAM \cite{TEAM} which also covers related technologies and areas. We shall not repeat this discussion here. Rather, we shall focus on the use of stationary vehicles to provision services, as work on this topic is less well-known but relevant for the research reported in this paper \cite{ibmICCVE}.

Information collected from stationary vehicles (i.e., vehicles with time-invariant position) has been used to study the performance of VANETs. For example, information from stationary vehicles has been used to accurately determine the stop-delay or idling times at road junctions. This information is then used to enhance road safety by detecting stop-line violations at signalised/controlled intersections or stopped vehicles inside tunnels, or to determine the availability of free parking spaces \cite{Mathur2010}. In the past, stationary vehicles have been mostly considered as passive nodes, and the relevant information was collected by using passive techniques such as proximity detectors and cameras. However, thanks to the advent of modern ITS technologies allowing vehicle-to-vehicle/infrastructure (V2X) communication and cooperative awareness, the use of stationary vehicles as active nodes to enhance services for VANETs is becoming a topic of great relevance. Some of the aforementioned services include: the improvement of multi-channel operations \cite{Campolo2012,Campolo2013}, the use as relays in content downloading and distribution \cite{Malandrino2012,Liu2012}, and mitigation of signal attenuation in ITS applications \cite{Sommer2013}.

Two main approaches can be distinguished when using stationary vehicles: the use of (i) powered-on vehicles, such as vehicles stopped in a queue; and (ii) the use of powered-off vehicles, such as parked vehicles. While the usage of powered-on stationary vehicles does not pose any particular technical requirement in the provision of new services, an obvious concern associated with parked vehicles is battery discharge; we shall have more to say about this shortly. An application highlighting the use of powered-on stationary vehicles is presented in \cite{Campolo2012} and \cite{Campolo2013}. In \cite{Sommer2013}, parked cars are postulated as relays in a multi-hop beaconing approach: they do not transmit their own information, but only retransmit information from moving vehicles. Battery consumption concerns are investigated in \cite{Crepaldi2013}, where the authors present a study on the impact of the communication system and processing unit on energy consumption. The main conclusion is that the power demanded from the communication system and the processing unit is not particularly critical. Specifically, given a fully charged battery, services can be provided for a significant time before a critical point of the charge level is reached (we shall have more to say on this in the next section).

%¶¶¶¶¶¶¶¶¶¶¶¶¶¶¶¶¶¶¶¶¶¶¶¶¶¶¶¶¶¶¶¶¶¶¶¶¶¶¶¶¶¶¶¶¶¶¶¶¶¶¶¶¶¶¶¶¶¶¶¶¶¶¶¶¶¶¶¶¶¶¶¶¶¶¶¶¶¶¶¶¶¶¶¶¶¶¶¶¶¶¶¶¶¶¶¶¶¶¶¶¶¶¶¶¶¶¶¶¶¶¶¶¶¶¶¶¶¶¶¶¶¶¶¶¶¶¶¶¶¶¶¶¶¶¶¶¶¶¶¶¶¶¶¶¶¶¶¶¶¶¶¶
%¶¶¶¶¶¶¶¶¶¶¶¶¶¶¶¶¶¶¶¶¶¶¶¶¶¶¶¶¶¶¶¶¶¶¶¶¶¶¶¶¶¶¶¶¶¶¶¶¶¶¶¶¶¶¶¶¶¶¶¶¶¶¶¶¶¶¶¶¶¶¶¶¶¶¶¶¶¶¶¶¶¶¶¶¶¶¶¶¶¶¶¶¶¶¶¶¶¶¶¶¶¶¶¶¶¶¶¶¶¶¶¶¶¶¶¶¶¶¶¶¶¶¶¶¶¶¶¶¶¶¶¶¶¶¶¶¶¶¶¶¶¶¶¶¶¶¶¶¶¶¶¶
\section{Statistics of stationary vehicles to support CP applications} \label{S:Characteristics}
%¶¶¶¶¶¶¶¶¶¶¶¶¶¶¶¶¶¶¶¶¶¶¶¶¶¶¶¶¶¶¶¶¶¶¶¶¶¶¶¶¶¶¶¶¶¶¶¶¶¶¶¶¶¶¶¶¶¶¶¶¶¶¶¶¶¶¶¶¶¶¶¶¶¶¶¶¶¶¶¶¶¶¶¶¶¶¶¶¶¶¶¶¶¶¶¶¶¶¶¶¶¶¶¶¶¶¶¶¶¶¶¶¶¶¶¶¶¶¶¶¶¶¶¶¶¶¶¶¶¶¶¶¶¶¶¶¶¶¶¶¶¶¶¶¶¶¶¶¶¶¶¶
%¶¶¶¶¶¶¶¶¶¶¶¶¶¶¶¶¶¶¶¶¶¶¶¶¶¶¶¶¶¶¶¶¶¶¶¶¶¶¶¶¶¶¶¶¶¶¶¶¶¶¶¶¶¶¶¶¶¶¶¶¶¶¶¶¶¶¶¶¶¶¶¶¶¶¶¶¶¶¶¶¶¶¶¶¶¶¶¶¶¶¶¶¶¶¶¶¶¶¶¶¶¶¶¶¶¶¶¶¶¶¶¶¶¶¶¶¶¶¶¶¶¶¶¶¶¶¶¶¶¶¶¶¶¶¶¶¶¶¶¶¶¶¶¶¶¶¶¶¶¶¶¶

A stationary vehicle without any hard time constraint for localisation, with a V2X communication system, and with access to Augmented GNSS (A-GNSS) based position estimation such as DGPS/AGPS \cite{Chan2012}, can easily fix its location precisely (with up to 10 cm of accuracy for Post-process DGPS). Thus it can become an anchor (anchor node) and, by broadcasting its location, it can help vehicles whose position is not known precisely (blind nodes) to localise themselves. Vehicles without hard time constraints include vehicles that are stationary for a long enough time, for possibly different reasons (e.g., traffic-light queues, bottlenecks, traffic jams and parked vehicles). We now give some basic statistics and facts to support the potential of stationary vehicles for CP applications.\newline

\textbf{1) Battery issues:} In the case of parked vehicles, the main technical challenge is keeping the on-board communication systems switched on when the vehicle is powered off. However, recent studies have shown that this is not a critical issue. In particular, by accessing a maximum of 10\% of the charge of a 480 Wh car's battery, it has been recently demonstrated that a 1W on-board unit can be constantly powered for up to 2 days \cite{Sommer2013}.\newline

\textbf{2) Availability of stationary powered-on vehicles:} The average time spent in stop/idle mode between two consecutive intersections in urban scenarios can reach up to 50\% of the whole time spent along the arterial road \cite{Lin2012}. Idle times become even larger at pm rush hours as people drive in a more relaxed way since they feel that they do not have strong time constraints in returning home, which causes longer queues at signalised/controlled intersections \cite{Sando2009}. Accordingly, these vehicles can be used with more priority for positioning applications.\newline

\begin{figure}[h]
\begin{centering}
\includegraphics[width=3.0in]{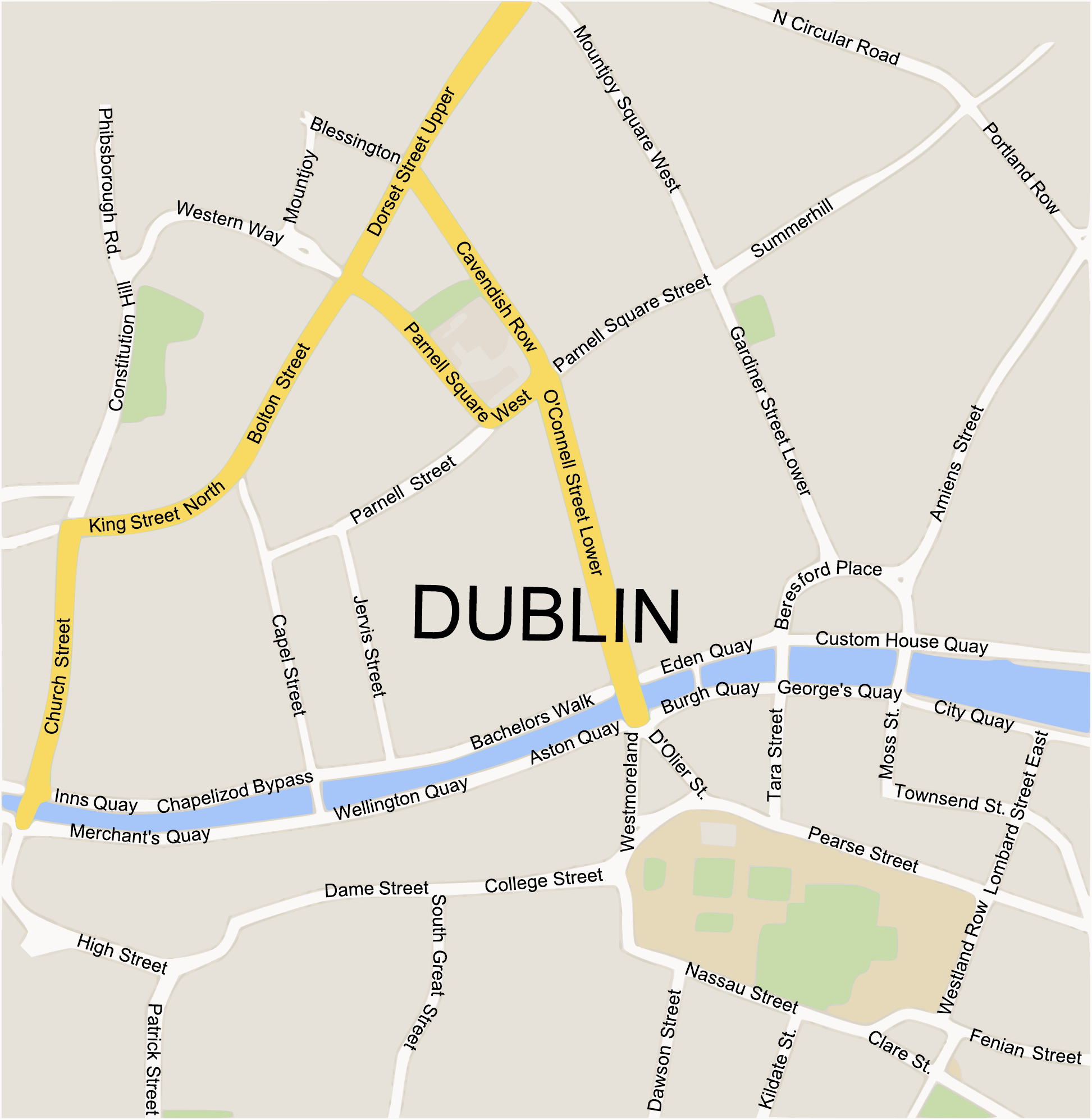}
\par\end{centering}
\caption{Roads in Dublin city center monitored between 16th-23rd June, 2014 using AA Roadwatch}.
\label{fig:Roads_Dublin}
\end{figure}

\textbf{3) Estimation of area occupancy of stationary powered-on vehicles:} In order to evaluate the potential impact of stationary powered-on vehicles on CP applications, it is important to evaluate the mobility area that is actually occupied by such vehicles with respect to the overall mobility network. As a support to our analysis, we use the statistics of the vehicular flows that are stuck or very slowly moving in the traffic (i.e., stationary powered-on vehicles), collected from the area of Dublin centre, Ireland, from AA Roadwatch (a service offered by AA Ireland\footnote{\url{http://www.theaa.ie/AA/AA-Roadwatch.aspx}}) as shown in Fig. \ref{fig:Roads_Dublin}.

\begin{figure}[h]
\begin{centering}
\includegraphics[width=3.0in]{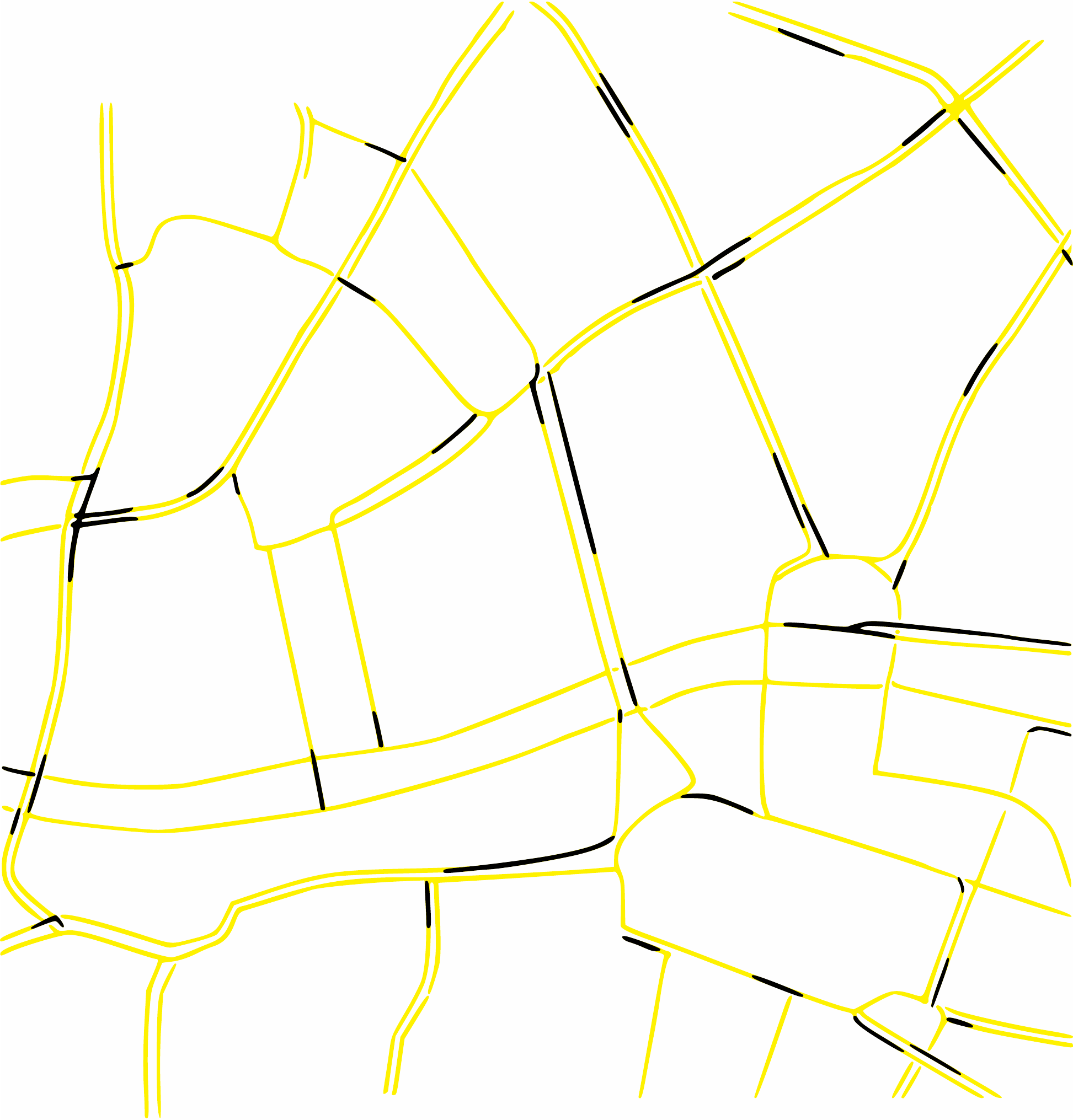}
\par\end{centering}
\caption{Measured data on June 17th, 2014, at the evening rush hour. Black sections: occupied by stationary powered-one vehicles; yellow sections: free traffic.}
\label{fig:Roads_Dublin_Tuesday}
\end{figure}

Data were collected two times per day between 16th-23rd June, 2014: one measure at the lunchtime rush hour (at 1 pm), and one at the evening rush hour (at 7:30 pm); the results are reported in Table \ref{tab:Results_AARoad}. From Table \ref{tab:Results_AARoad} it is possible to conclude that, in the context of this small survey, an average of 22.64\% and 15.71\% of the main streets in Dublin's city centre are occupied by stationary powered-on vehicles at the rush hours of working days. Even though numbers decrease in the weekends to 14.36\% (at lunchtime), and to 13.03\% (in the evening), such numbers suggest the benefit of providing support for the localisation of blind vehicles during rush hours. In addition, note that the covered area is concentrated around intersections, as shown with black curves in Fig. \ref{fig:Roads_Dublin_Tuesday}, where most of the cars (and thus, of the blind nodes as well) are also concentrated.\newline

\begin{table}[h]
\caption{Percentage of slow and very slow moving traffic overall traffic for the main streets in Dublin city center as shown in Figure \ref{fig:Roads_Dublin}.}
\label{tab:Results_AARoad}
\begin{centering}
\begin{tabular}{|c|c|c|}
\cline{2-3}
\multicolumn{1}{c|}{} & \multicolumn{2}{c|}{Rush hour}\tabularnewline
\hline
Day & Lunchtime & Evening\tabularnewline
\hline
\hline
Monday & 18.39\% & 13.07\% \tabularnewline
\hline
Tuesday & 17.68\% & 16.17\% \tabularnewline
\hline
Wednesday & 25.03\% & 14.58\% \tabularnewline
\hline
Thursday & 23.49\% & 20.67\% \tabularnewline
\hline
Friday & 28.63\% & 14.07\% \tabularnewline
\hline
Saturday & 14.95\% & 9.72\% \tabularnewline
\hline
Sunday & 13.77\% & 16.33\% \tabularnewline
\hline
\end{tabular}
\par\end{centering}
\end{table}

\textbf{4) Availability of stationary powered-off vehicles:} A car is typically parked on average up to 23 hours a day \cite{Pasaoglu2012,Litman2013}, and usually outdoor \cite{Malandrino2012,Sommer2013}. For example, in a recent study of 61,000 daily parking events in Montreal City \cite{Morency2008}, 69.2\% of all parked cars were parked on-street and 27.1\% were parked in outside parking lots, while only 3.7\% were parked in interior parking facilities, and the average duration of on-street parking was 6.64 hours. Further, current predictions concerning regulated parking spaces claim that the average portion of on-street parking spaces, as a percentage of the overall total number of parking spaces, will be 30.17\% (with up to 56.22\% for Italy and 43.30\% in Spain) \cite{EPA2013}. Thus ``parked cars'' can be thought of as a (so far) unused and dense sensor network with no power or collection constraints, and whose position can be known very precisely. Consequently, such a network of static nodes can be used to help in solving the CP problem for ITS systems.\newline

\begin{table}[h]
\begin{centering}
\caption{The United States Federal Communication Commission (FCC) classification for DSRC devices \cite{Kenney2011}.}
\label{tab:Table_DSRC}
\begin{tabular}{|c|c|c|}
\hline
Device & Max. Output & Communication\tabularnewline
class & Power (dBm) & zone (meters)\tabularnewline
\hline
\hline
A & 0 & 15\tabularnewline
\hline
B & 10 & 100\tabularnewline
\hline
C & 20 & 400\tabularnewline
\hline
D & 28.8 & 1000\tabularnewline
\hline
\end{tabular}
\par\end{centering}
\end{table}

\begin{figure}[h]
\centering
\subfloat[]{\includegraphics[width=1.7in]{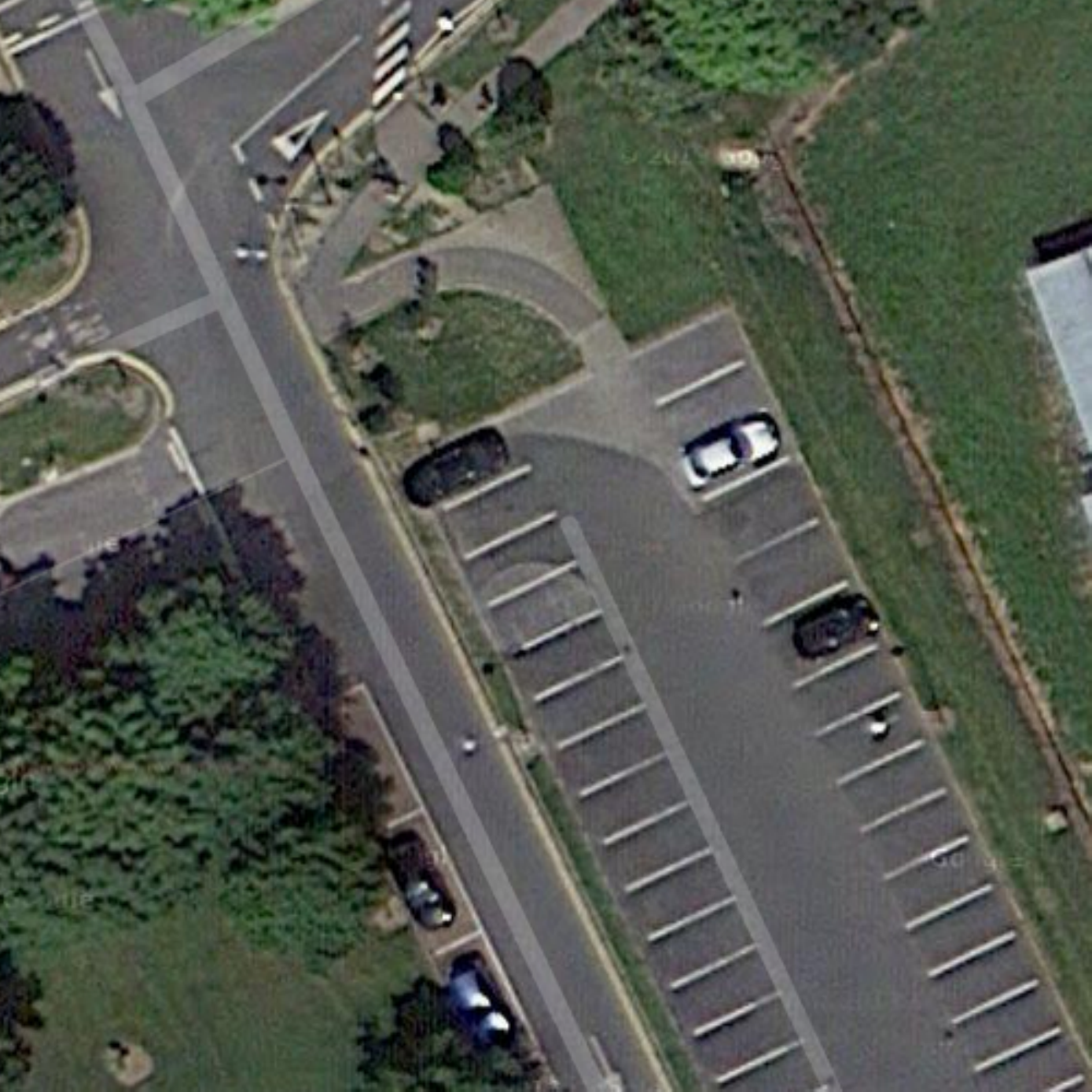}
\label{fig:Fig03a}}
\subfloat[]{\includegraphics[width=1.7in]{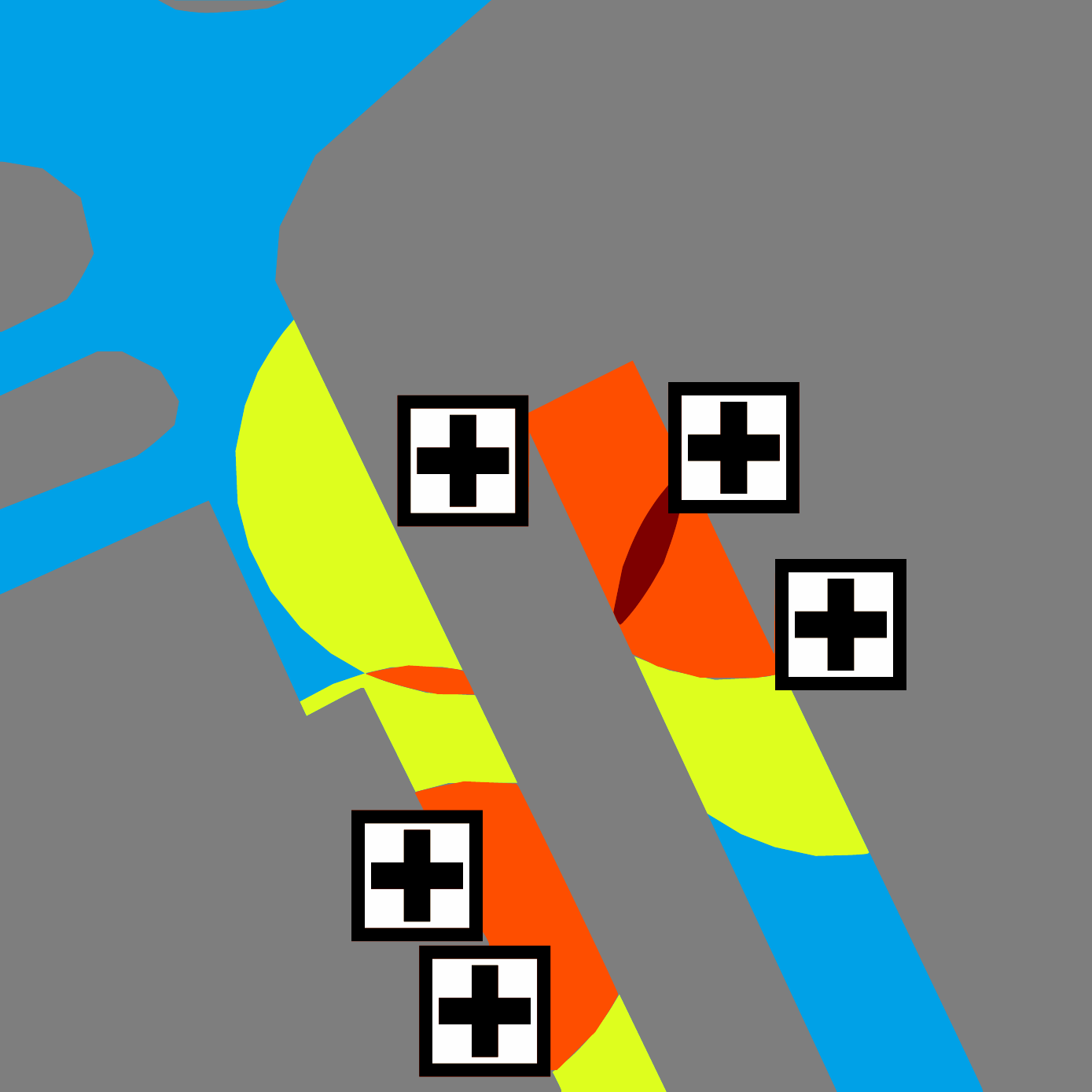}
\label{fig:Fig03b}}
\caption{Illustrative scenario for area coverage estimation: a) selected region, and b) identified areas: vehicular transit (blue), parked vehicles (black crosses),
and covered area for a radius equal to 15 m. Color convention for levels of coverage (see Table \ref{tab:Levels_Coverage}): yellow for Level 1, orange for Level 2, and red for Level 3.}
\label{fig:Scenario_small}
\end{figure}

\textbf{5) Estimation of the area covered by stationary powered-off vehicles:} The expected size of the area that can be covered by parked vehicles does clearly depend on the specific technology used to broadcast the required information, and on the amount of cars parked in the area.
Table \ref{tab:Table_DSRC} summarizes the characteristics of the most popular Dedicated Short-Range Communication (DSRC) technologies. For instance, we are going to estimate the area coverage as follows:
~
\begin{enumerate}
\item Identify area of vehicular transit (blue area in Fig. \ref{fig:Scenario_small}.b), referred as $A_{Transit}$;
\item Identify stationary vehicles (black crosses in Fig. \ref{fig:Scenario_small}.b) and their total circular coverage (coloured circles in Fig. \ref{fig:Scenario_small}.b), referred as $A_{Coverage}$;
\item Calculate the intersection between the $A_{Coverage}$ and $A_{Transit}$, referred to as $A_{C\cap T}$.
\item Calculate the ratio between $A_{C\cap T}$ and $A_{Transit}$ as a function of the levels given in Table \ref{tab:Levels_Coverage}.
\end{enumerate}

\begin{table}[h]
\begin{centering}
\caption{Proposed levels of coverage.}
\label{tab:Levels_Coverage}
\begin{tabular}{|c|c|c|}
\hline
Level of coverage & Description \tabularnewline
\hline
\hline
Level 1 & Coverage with only 1 stationary car \tabularnewline
\hline
Level 2 & Coverage with only 2 stationary cars \tabularnewline
\hline
Level 3 & Coverage with 3 or more stationary cars \tabularnewline
\hline
\end{tabular}
\par\end{centering}
\end{table}

\begin{figure}[h]
\begin{centering}
\includegraphics[width=2.5in]{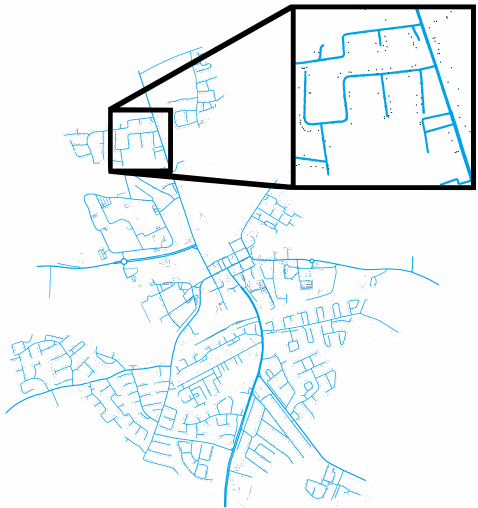}
\par\end{centering}
\caption{Allowed area for vehicular transit (blue) and parked vehicles (black) in Maynooth, Co. Kildare, Dublin, using satellite imagery from Google Maps. A total of 4541 vehicles were identified as parked vehicles.}
\label{fig:Road_parked}
\end{figure}

\begin{figure}[h]
\begin{centering}
\includegraphics[width=2.5in]{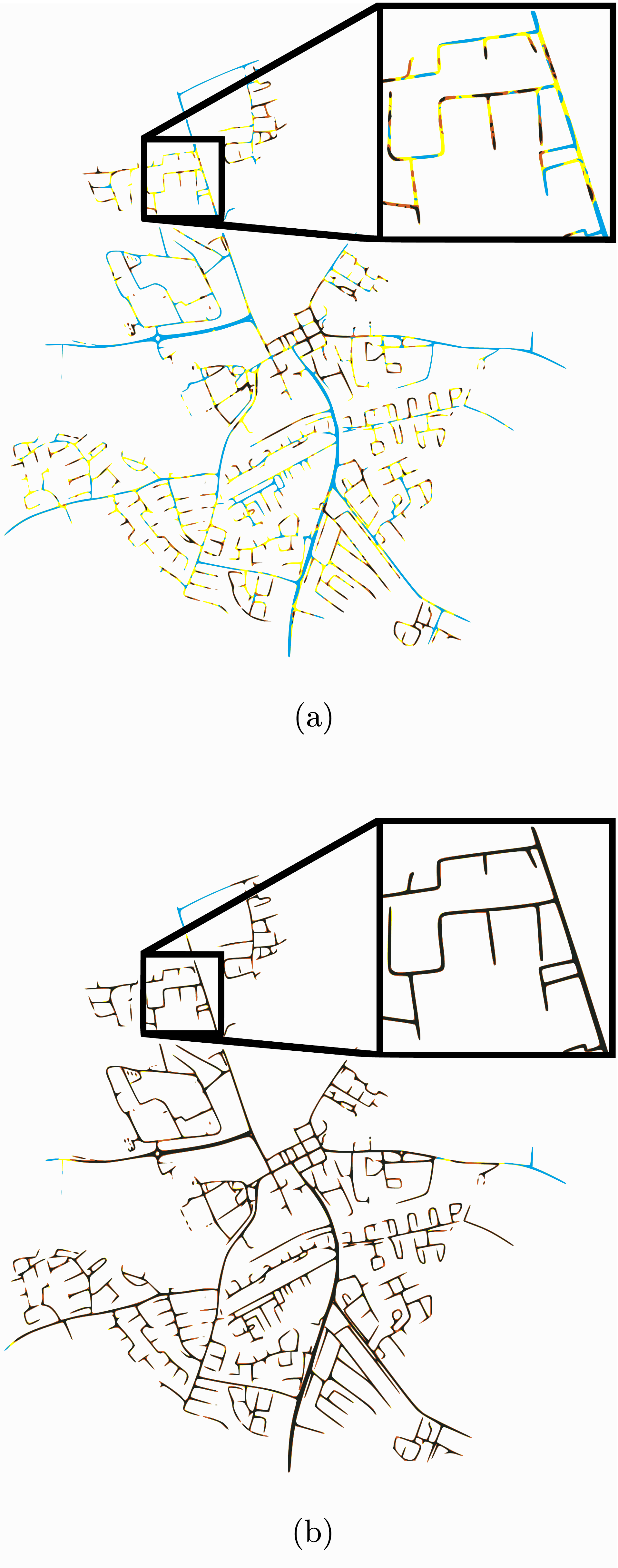}
\par\end{centering}
\caption{Covered areas for the case study using different communication zones: a) 15 m (Class-A DSRC devices), and b) 100 m (Class-B DSRC devices). Color convention for levels of coverage (see Table \ref{tab:Levels_Coverage}): yellow for Level 1, orange for Level 2, and red for Level 3.}
\label{fig:Corevage_015_100}
\end{figure}

For illustration of this approach, we used top-view images from Google Maps of Maynooth Town (co. Kildare, Ireland) with an image from October 31st 2013.
The parked vehicles and the area coverage using the above procedure are shown in Fig. \ref{fig:Road_parked} and Fig. \ref{fig:Corevage_015_100} respectively, and the coverage percentage is summarised in Table \ref{tab:Areas}.

\begin{table}[h]
\caption{Potential covered areas for the case study. The 100\% corresponds to the areas allowed for vehicular transit.}
\label{tab:Areas}
\begin{centering}
\begin{tabular}{|l|c|c|}
\cline{2-3}
\multicolumn{1}{l|}{} & \multicolumn{2}{c|}{Com. zone (meters)}\tabularnewline
\cline{2-3}
\multicolumn{1}{l|}{} & 15 & 100 \tabularnewline
\hline
Covered area of Level 3& 27.78\% & 97.02\%\tabularnewline
\hline
Covered area of Level 2 & 11.98\% & 0.22\%\tabularnewline
\hline
Covered area of Level 1& 14.81\% & 0.31\%\tabularnewline
\hline
No covered area & 45.43\% & 2.45\%\tabularnewline
\hline
\end{tabular}
\par\end{centering}
\end{table}

In this simple example, according to Table \ref{tab:Areas}, $27.78\%$ of the allowed area for vehicular transit can be covered using Class-A DSRC devices, and $97.02\%$ using Class-B DSRC devices. Thus, parked vehicles can potentially cover most of the travelling area.\\
\newline
\textbf{Comment 1: } Note that in this study we are not explicitly considering multipath effects and possible shadowing consequences.
Rather, we assume a direct line-of-sight between moving and stationary vehicles. In such circumstances, such effects are mitigated \cite{GiovanniBellusci2009,Gozalvez2012}. From a
practical perspective, this means in the remainder of the paper, we assume mainly short-range communications between vehicles with a communication range of 15 m - the
underlying assumption being direct line-of-sight. As we see from Table \ref{tab:Areas}, even with this assumption, the total covered area in our experiments is still significant (more than $54\%$).\\

Finally, note that the total coverage is actually greater, once one considers the stationary powered-on vehicles as well. Since, as already remarked, stationary powered-on vehicles can mainly contribute to cover the intersection areas, while stationary powered-off vehicles are usually densely distributed along the roads, in principle a full coverage of the whole urban mobility network can be potentially achieved.

%¶¶¶¶¶¶¶¶¶¶¶¶¶¶¶¶¶¶¶¶¶¶¶¶¶¶¶¶¶¶¶¶¶¶¶¶¶¶¶¶¶¶¶¶¶¶¶¶¶¶¶¶¶¶¶¶¶¶¶¶¶¶¶¶¶¶¶¶¶¶¶¶¶¶¶¶¶¶¶¶¶¶¶¶¶¶¶¶¶¶¶¶¶¶¶¶¶¶¶¶¶¶¶¶¶¶¶¶¶¶¶¶¶¶¶¶¶¶¶¶¶¶¶¶¶¶¶¶¶¶¶¶¶¶¶¶¶¶¶¶¶¶¶¶¶¶¶¶¶¶¶¶
%¶¶¶¶¶¶¶¶¶¶¶¶¶¶¶¶¶¶¶¶¶¶¶¶¶¶¶¶¶¶¶¶¶¶¶¶¶¶¶¶¶¶¶¶¶¶¶¶¶¶¶¶¶¶¶¶¶¶¶¶¶¶¶¶¶¶¶¶¶¶¶¶¶¶¶¶¶¶¶¶¶¶¶¶¶¶¶¶¶¶¶¶¶¶¶¶¶¶¶¶¶¶¶¶¶¶¶¶¶¶¶¶¶¶¶¶¶¶¶¶¶¶¶¶¶¶¶¶¶¶¶¶¶¶¶¶¶¶¶¶¶¶¶¶¶¶¶¶¶¶¶¶
\section{Embedding stationary vehicles into a CP technique} \label{S:CP_approach}
%¶¶¶¶¶¶¶¶¶¶¶¶¶¶¶¶¶¶¶¶¶¶¶¶¶¶¶¶¶¶¶¶¶¶¶¶¶¶¶¶¶¶¶¶¶¶¶¶¶¶¶¶¶¶¶¶¶¶¶¶¶¶¶¶¶¶¶¶¶¶¶¶¶¶¶¶¶¶¶¶¶¶¶¶¶¶¶¶¶¶¶¶¶¶¶¶¶¶¶¶¶¶¶¶¶¶¶¶¶¶¶¶¶¶¶¶¶¶¶¶¶¶¶¶¶¶¶¶¶¶¶¶¶¶¶¶¶¶¶¶¶¶¶¶¶¶¶¶¶¶¶¶
%¶¶¶¶¶¶¶¶¶¶¶¶¶¶¶¶¶¶¶¶¶¶¶¶¶¶¶¶¶¶¶¶¶¶¶¶¶¶¶¶¶¶¶¶¶¶¶¶¶¶¶¶¶¶¶¶¶¶¶¶¶¶¶¶¶¶¶¶¶¶¶¶¶¶¶¶¶¶¶¶¶¶¶¶¶¶¶¶¶¶¶¶¶¶¶¶¶¶¶¶¶¶¶¶¶¶¶¶¶¶¶¶¶¶¶¶¶¶¶¶¶¶¶¶¶¶¶¶¶¶¶¶¶¶¶¶¶¶¶¶¶¶¶¶¶¶¶¶¶¶¶¶

The objective of this section is to illustrate how any available technique to realise cooperative positioning (e.g., triangulation, trilateration, cooperative estimation) can be improved by further using stationary vehicles.

%============================================================================
\subsection{Using a stationary vehicle as an anchor node}
%============================================================================

A stationary vehicle can act as an anchor node to support the localisation of other vehicles, only after its own position is precisely known. Note that powered-off stationary vehicles (e.g., parked vehicles) do not have hard time constraints for localisation, and they can localise themselves in an accurate manner using any (possibly slow) localisation technique. For instance, they can localise themselves using A-GNSS-based position estimation and V2I communication systems, or if they have access to enough information from some nearby anchor nodes via V2V communication, then multilateration can be used for localisation. Once a stationary vehicle has a precise fix on its location, then it can be used as an anchor node to support the CP of other vehicles. Assuming that some stationary vehicles exist, we shall now provide algorithms to localise other vehicles using the stationary vehicles as anchor nodes.

%============================================================================
\subsection{Localising blind stationary vehicles}\label{s:Localising_blind_stationary}
%============================================================================

A blind stationary vehicle is a stationary vehicle whose localisation is not yet precise. To localise blind stationary vehicles, the following strategy is adopted:
\begin{itemize}
\item If 3 or more anchor neighbours are available:
\begin{itemize}
\item (\underline{default option}): localise using a CP technique supported by V2V communication to anchor nodes;
\item (\underline{back-up option}): localise over longer time horizon using blue A-GNSS-based position estimation.
\end{itemize}

\item If 2 anchor neighbours are available:
\begin{itemize}
\item (\underline{default option}): localise using a CP technique supported by V2V communication to anchor nodes and additional information
from other surrounding vehicles (e.g., previous positions/speeds/ranges) to solve ambiguities;
\item (\underline{back-up option}): localise over longer time horizon using A-GNSS-based position estimation.
\end{itemize}

\item Otherwise:
\begin{itemize}
\item (\underline{default option}): localise only using A-GNSS-based position estimation.
\end{itemize}

\end{itemize}

In addition, a blind powered-off stationary vehicle remains as an inactive node in the VANET always before becoming an anchor. Also recall that a blind stationary
vehicle can only become an anchor after having reached a high accuracy in the localisation procedure (e.g., when having access to enough
anchors, or via A-GNSS-based position estimation).\\
\newline
\textbf{Comment 2: }{The interested reader will notice that some blind nodes will eventually know their position in a more accurate fashion than others, mainly depending on both the quality of the inter-vehicle distance measurements and the specific localisation procedure. However, in the sequel, we shall ignore this transient effect for powered-off stationary cars, and assume that all remaining stationary cars localise themselves accurately and quickly. Only precisely located cars are used as anchors; the effects of different levels of accuracy of anchor nodes is beyond the scope of this present paper.}

%============================================================================
\subsection{Localising blind moving vehicles}
%============================================================================

The procedure is similar to that given before. However, in this case, at the end of the localisation process, moving vehicles are not considered as anchor nodes. A blind moving vehicle with access to accurate information from anchor nodes, will become a {\it pseudo-anchor} node, and will be used with certain priority in the localisation of other vehicles, if needed.

%........................................................................................................................................................
\subsection{Node selection strategy for localisation}
%........................................................................................................................................................

A smart selection of the neighbouring nodes for CP can speed up the localisation procedure and mitigate energy requirements, and thus, battery consumption. In practical cases, a large number of anchor nodes can be used for localisation, and strategies to select the ``best'' subset of nodes have been already investigated in the literature. See for example \cite{Zhang2012b} where an optimal subset of anchor nodes is identified according to a geometric dilution of precision (GDOP) process. Here, we adopt
a simple deterministic selection strategy that assigns priorities to anchor nodes, pseudo-anchor nodes and blind nodes respectively. Although in principle any other method for node selection can be used, we use the following procedure.
\begin{itemize}
\item Any surrounding car $c_i$ has a priority level $p$ given by
\begin{equation}
p\left(c_i\right)=\left\{ \begin{array}{cc}
1, & \mbox{if \ensuremath{c_i} is an anchor node,}\\
2, & \mbox{if \ensuremath{c_i} is a pseudo-anchor node,}\\
3, & \mbox{otherwise (i.e., blind nodes).}
\end{array}\right.\label{eq:Priority}
\end{equation}
where 1 is the highest priority.
\end{itemize}

Three surrounding vehicles inside the communication zone of the target vehicle are sufficient for the localisation process. Accordingly, the selection strategy has to choose the three ``best'' heighbours, for which we recommend those with higher priority should be selected first. Given the availability of vehicles with the same priority, then we choose those that are closest to the vehicle to be localised. If fewer than three are available, we choose the highest number of available cars.

%........................................................................................................................................................
\subsection{Integrating stationary vehicles into any CP algorithm}
%........................................................................................................................................................

Stationary vehicles have the potential to improve CP algorithms in several ways:
\begin{itemize}
\item[(i)] Parked vehicles augment fixed roadside-like infrastructure, and so more accurate and reliable range information is available to moving vehicles. Thus, parked vehicles have the potential to greatly augment {\it range based-}, {\it fine-grained- }, {\it anchor-based- }, {\it absolute- }, {\it sequential-}, and {\it distributed CP algorithms}.
\item[(ii)] Incorporating parked vehicles into CP algorithms implies that even when other moving vehicles or dedicated infrastructure are not available for positioning, cooperative positioning may still be possible by using parked vehicles. This greatly helps a host of CP algorithms, including {\it range-free-}, {\it coarse-grained-}, {\it anchor-free-}, {\it relative-}, {\it multi-hop-}, {\it concurrent-}, and {\it centralised CP algorithms} (the latter in small-scale networks).
\item[(iii)] As parked vehicles are usually in close proximity to the road, there might be an advantage with respect to noise properties over dedicated infrastructure for: {\it single-hop-}, {\it line-of-sight-}, {\it Bayesian CP algorithms}, and {\it standard/modern estimation algorithms for CP}.
\item[(iv)] Finally, localised stationary vehicles can help to localise other stationary vehicles with no position information. Thus, they augment {\it CP algorithms for static networks}.
\end{itemize}

%¶¶¶¶¶¶¶¶¶¶¶¶¶¶¶¶¶¶¶¶¶¶¶¶¶¶¶¶¶¶¶¶¶¶¶¶¶¶¶¶¶¶¶¶¶¶¶¶¶¶¶¶¶¶¶¶¶¶¶¶¶¶¶¶¶¶¶¶¶¶¶¶¶¶¶¶¶¶¶¶¶¶¶¶¶¶¶¶¶¶¶¶¶¶¶¶¶¶¶¶¶¶¶¶¶¶¶¶¶¶¶¶¶¶¶¶¶¶¶¶¶¶¶¶¶¶¶¶¶¶¶¶¶¶¶¶¶¶¶¶¶¶¶¶¶¶¶¶¶¶¶¶
%¶¶¶¶¶¶¶¶¶¶¶¶¶¶¶¶¶¶¶¶¶¶¶¶¶¶¶¶¶¶¶¶¶¶¶¶¶¶¶¶¶¶¶¶¶¶¶¶¶¶¶¶¶¶¶¶¶¶¶¶¶¶¶¶¶¶¶¶¶¶¶¶¶¶¶¶¶¶¶¶¶¶¶¶¶¶¶¶¶¶¶¶¶¶¶¶¶¶¶¶¶¶¶¶¶¶¶¶¶¶¶¶¶¶¶¶¶¶¶¶¶¶¶¶¶¶¶¶¶¶¶¶¶¶¶¶¶¶¶¶¶¶¶¶¶¶¶¶¶¶¶¶
\section{Performance evaluation} \label{S:Evaluation}
%¶¶¶¶¶¶¶¶¶¶¶¶¶¶¶¶¶¶¶¶¶¶¶¶¶¶¶¶¶¶¶¶¶¶¶¶¶¶¶¶¶¶¶¶¶¶¶¶¶¶¶¶¶¶¶¶¶¶¶¶¶¶¶¶¶¶¶¶¶¶¶¶¶¶¶¶¶¶¶¶¶¶¶¶¶¶¶¶¶¶¶¶¶¶¶¶¶¶¶¶¶¶¶¶¶¶¶¶¶¶¶¶¶¶¶¶¶¶¶¶¶¶¶¶¶¶¶¶¶¶¶¶¶¶¶¶¶¶¶¶¶¶¶¶¶¶¶¶¶¶¶¶
%¶¶¶¶¶¶¶¶¶¶¶¶¶¶¶¶¶¶¶¶¶¶¶¶¶¶¶¶¶¶¶¶¶¶¶¶¶¶¶¶¶¶¶¶¶¶¶¶¶¶¶¶¶¶¶¶¶¶¶¶¶¶¶¶¶¶¶¶¶¶¶¶¶¶¶¶¶¶¶¶¶¶¶¶¶¶¶¶¶¶¶¶¶¶¶¶¶¶¶¶¶¶¶¶¶¶¶¶¶¶¶¶¶¶¶¶¶¶¶¶¶¶¶¶¶¶¶¶¶¶¶¶¶¶¶¶¶¶¶¶¶¶¶¶¶¶¶¶¶¶¶¶

In order to evaluate the proposed approach, we used simulations in SUMO, a popular mobility simulator\footnote{\url{www.dlr.de/ts/sumo/en/}, commonly used in the ITS community}. In the following, we consider both stationary powered-on and powered-off cars represented by not moving cars in traffic and parked cars, respectively. Moreover, DGPS/AGPS are assumed to be running in parallel on each vehicle (see Section \ref{s:Localising_blind_stationary}), and we also assume that all parked vehicles are available to be used as anchor nodes from the beginning of all simulations.

\subsection{CP algorithms to be evaluated}

The first CP algorithm to be evaluated is a decentralised estimation algorithm based on Guaranteed Convergence Particle Swarm Optimiser (GCPSO) \cite{VandenBergh2010},
using inter-vehicle distance measurements to minimise the cost function
\[
f\left(\hat{\mathbf{x}}_{i,t}\right)=\sum_{\forall j\in S}\left\Vert d_{ij,t}-\left\Vert \mathbf{\hat{x}}_{i,t}-\mathbf{x}_{j,t}\right\Vert \right\Vert ^{2}+\left\Vert \mathbf{\hat{x}}_{i,t}-T_{s}\mathbf{v}_{i,t-1}\right\Vert ^{2},
\]
where $\hat{\mathbf{x}}_{i,t}$ is the estimate of the real (unavailable) 2D position ${\mathbf{x}}_{i}$ of car $i$ at time step $t$,
$S$ is the set of all one-hop neighbour cars of car $i$ inside its communication zone, $\mathbf{x}_j$ is the available 2D position for car $j$ shared to car $i$,
${\mathbf{v}}_{i}$ is the 2D speed of car $i$, $T_s$ is the sampling time for data collection,
and $d_{ij}$ is the measured distance between cars $i$ and $j$,
where $j$ is the index of a car in $S$. Thus, this is a distributed, range-based CP algorithm based on a heuristic optimiser.

The second CP algorithm to be evaluated is the conventional Extended Kalman Filter (EKF) \cite{Parker2007}. The implementation in \cite{Parker2007}
fuses inter-vehicle distance measurements and vehicle kinematics (speed information), to obtain a sequential Bayesian estimation of the cars' position.
Thus, this is a distributed, range-based, Bayesian CP algorithm. For a detailed description of this algorithm please refer to \cite{Parker2007}.

In our simulations we use the following parameter settings and make the following assumptions:
\begin{itemize}
\item General assumptions \cite{Parker2007,GiovanniBellusci2009}:

\begin{itemize}
\item sampling time $T_s=1s$, ranging-error and GPS noise assumed as zero-mean Gaussian variable with standard deviations $\sigma_{\mathbf{R}}\in\{0.2,4\}$
and $\sigma_{GPS} = 6$, respectively;
\end{itemize}

\item GCPSO \cite{VandenBergh2010}:
\begin{itemize}
\item number of particles = 4, number of iterations = 20, $s_c = 15$, $f_c = 5$, $\rho = 1$,
$c_{1,2}=2$, linearly decreasing inertia weight from 0.9 to 0.2, minimum fitness value for stopping criteria = 0,
initial swarm: 2 particles at the position of the neighbour with the highest priority and
2 particles at $T_{s}\mathbf{v}_{i,t-1}$;
\end{itemize}

\item EKF \cite{Parker2007}:
\begin{itemize}
\item $\mathbf{Q}_{k-1}=\sigma_{\mathbf{Q}}^2\mathbf{I}$, $\mathbf{R}_{k}=\sigma_{\mathbf{R}}^2\mathbf{I}$, $\mathbf{\Gamma}_{k-1}=\sigma_{\mathbf{\Gamma}}^2\mathbf{I}$, with $\mathbf{I}$
the identity matrix with proper dimensions, and with $\sigma_{\mathbf{\Gamma}}=0.5$, $\sigma_{\mathbf{Q}} = 2$ and $\sigma_{\mathbf{R}}$ as in general assumptions.
\end{itemize}

\end{itemize}

All the previous assumptions, and the choice of parameters, are consistent with those used in \cite{Parker2007,GiovanniBellusci2009,VandenBergh2010},
where $\sigma_{\mathbf{R}}=0.2$ corresponds to UBW-based ranging and $\sigma_{\mathbf{R}}=4$ corresponds to radio-based ranging.

Finally, the position estimation for isolated cars (i.e., those with no surrounding neighbours) is updated using their past position estimates and new velocity
measurements. Since this method leads to an increasing localisation error due to cumulative effects, we occasionally use the measured GPS position to avoid
long term effects from this kind of error in isolated cars.

%============================================================================
\subsection{Small scale scenario}\label{s:small_scale_example}
%============================================================================
\label{Small_Scale}

\begin{figure}[h]
\begin{centering}
\includegraphics[width=3.0in]{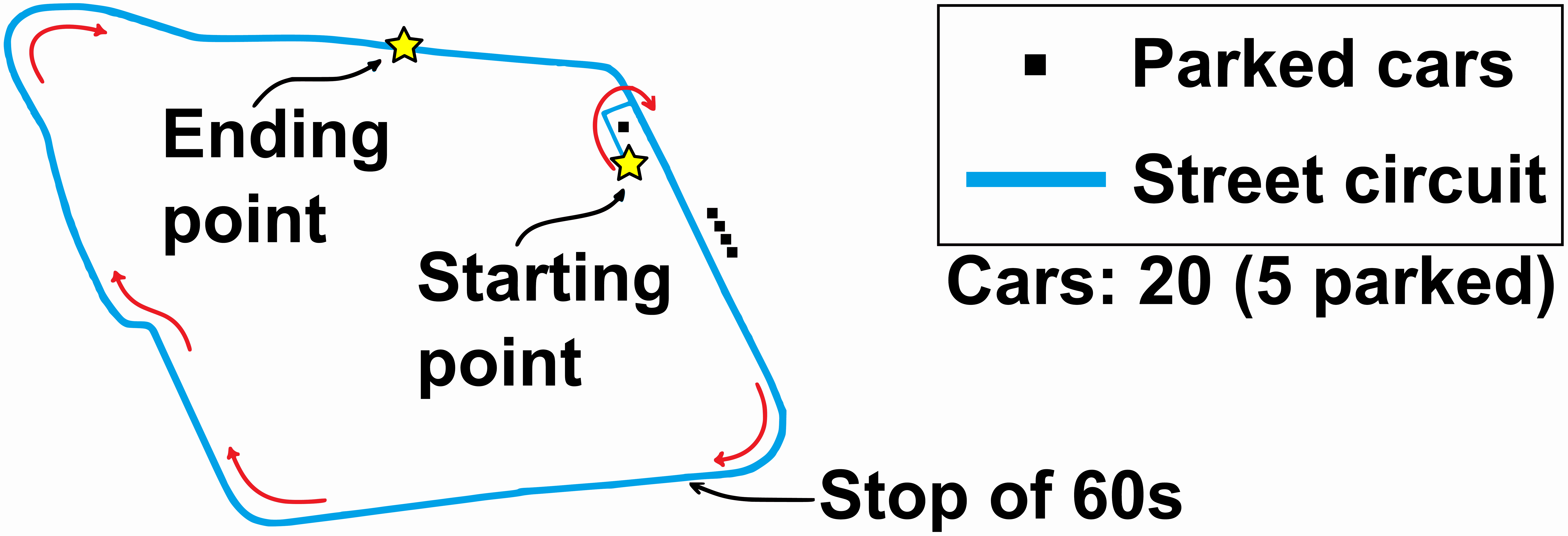}
\par\end{centering}
\caption{Circuit for the small scale test (blue) and selected parked cars (black).}
\label{fig:Street_circuit_blue}
\end{figure}

The first scenario to be tested is a street circuit inside the NUIM North Campus (Fig. \ref{fig:Street_circuit_blue}). Here, the main
idea is to evaluate the performance of the proposed approach in terms of the range of communication zones for V2V communication corresponding
to Class-A/B DSRC devices (see Table \ref{tab:Table_DSRC}), i.e., 15 m and 100 m respectively. The results obtained for an ensemble of 40 experiments per CP approach
are given in Table \ref{tab:RMSE_Example1}.

In Table \ref{tab:RMSE_Example1} we observe that an average improvement of up to 9.92 \% in the localisation RMSE
can be obtained when using stationary vehicles as prioritised nodes in the CP process and a communication zone of 15 m. This result allow us to
conclude that a communication zone of even 15 m potentially offers a significant improvement for the proposed CP approach over
the traditional approaches in particular, even when there are small numbers of parked vehicles.

Note also that with an increased communication range of 100 m, significant improvements can be obtained.
For example, this increased communication range allows, for our example,  an improvement\footnote{Concerning the proposed approach with respect to the traditional approach.}
 of up to $49.04\%$ in the localisation RMSE, with an accuracy level up to $3.26$~m (EKF with $\sigma_R=0.2$).
Clearly, in this case, multipath and shadowing effects cannot be neglected.

\begin{table}[h]
\caption{Results for the small scale scenario using two different values for the communication zone (CZ):
Localisation RMSE and improvements for the target car.}
\label{tab:RMSE_Example1}

\begin{centering}
\begin{tabular}{|l|c|c|c|c|c|c|c|}
\cline{4-7} 
\multicolumn{1}{l}{} & \multicolumn{1}{c}{} &  & \multicolumn{4}{c|}{RMSE (meters)} & \multicolumn{1}{c}{}\tabularnewline
\cline{4-7} 
\multicolumn{1}{l}{} & \multicolumn{1}{c}{} &  & \multicolumn{2}{c|}{Traditional CP} & \multicolumn{2}{c|}{Proposed CP} & \multicolumn{1}{c}{}\tabularnewline
\cline{3-8} 
\multicolumn{1}{l}{} &  & \multirow{2}{*}{$\sigma_{R}$} & \multirow{2}{*}{Mean} & \multirow{2}{*}{$\sigma$} & \multirow{2}{*}{Mean} & \multirow{2}{*}{$\sigma$} & Average\tabularnewline
\multicolumn{1}{l}{} &  &  &  &  &  &  & improvement\tabularnewline
\hline 
\hline 
\multirow{4}{*}{\begin{sideways}
CZ: 15 m
\end{sideways}} & \multirow{2}{*}{GCPSO} & 0.2  & 6.92 & 4.05 & 6.31 & 3.75 & 8.74\%\tabularnewline
\cline{3-8} 
 &  & 4  & 8.00 & 4.37 & 7.62 & 4.33 & 4.82\%\tabularnewline
\cline{2-8} 
 & \multirow{2}{*}{EKF} & 0.2 & 7.63 & 4.47 & 6.88 & 4.50 & 9.92\%\tabularnewline
\cline{3-8} 
 &  & 4 & 7.78 & 4.17 & 7.18 & 4.12 & 7.65\%\tabularnewline
\hline 
\hline 
\multirow{4}{*}{\begin{sideways}
CZ: 100 m
\end{sideways}} & \multirow{2}{*}{GCPSO} & 0.2  & 6.84 & 4.12 & 3.49 & 2.57 & 49.04\%\tabularnewline
\cline{3-8} 
 &  & 4  & 13.69 & 8.19 & 8.22 & 5.46 & 39.92\%\tabularnewline
\cline{2-8} 
 & \multirow{2}{*}{EKF} & 0.2  & 5.98 & 3.27 & 3.26 & 2.26 & 45.45\%\tabularnewline
\cline{3-8} 
 &  & 4  & 7.65 & 4.34 & 5.15 & 3.32 & 32.71\%\tabularnewline
\hline 
\end{tabular}
\par\end{centering}

~

CZ: Communication zone.

Mean: average of 40 different measurements.

$\sigma$: standard deviation.

\end{table}

%\\

%\begin{table}[h]
%\caption{Results for the small scale scenario using 100 m of communication zone:
%Localisation RMSE and improvements for the target car.}
%\label{tab:RMSE_Example1_100}
%
%\begin{centering}
%\begin{tabular}{|c|c|c|c|c|c|c|}
%\cline{3-6}
%\multicolumn{1}{c}{} & & \multicolumn{4}{c|}{RMSE (meters)} & \multicolumn{1}{c}{}\tabularnewline
%\cline{3-6}
%\multicolumn{1}{c}{} & & \multicolumn{2}{c|}{Traditional CP} & \multicolumn{2}{c|}{Proposed CP} & \multicolumn{1}{c}{}\tabularnewline
%\cline{3-7}
%\multicolumn{1}{c}{} & & \multirow{2}{*}{Mean} & \multirow{2}{*}{$\sigma$} & \multirow{2}{*}{Mean} & \multirow{2}{*}{$\sigma$} & Average\tabularnewline
%\cline{2-2}
%\multicolumn{1}{c|}{} & $\sigma_{R}$ & & & & & improvement\tabularnewline
%\hline
%\hline
%\multirow{2}{*}{GCPSO} & 0.2 & 6.84 & 4.12 & 3.49 & 2.57 & 49.04\%\tabularnewline
%\cline{2-7}
%& 4 & 13.69 & 8.19 & 8.22 & 5.46 & 39.92\%\tabularnewline
%\hline
%\hline
%\multirow{2}{*}{EKF} & 0.2 & 5.98 & 3.27 & 3.26 & 2.26 & 45.45\%\tabularnewline
%\cline{2-7}
%& 4 & 7.65 & 4.34 & 5.15 & 3.32 & 32.71\%\tabularnewline
%\hline
%\end{tabular}
%\par\end{centering}
%
%~
%
%Mean: average of 40 different measurements.
%
%$\sigma$: standard deviation.
%
%\end{table}

%============================================================================
\subsection{Large scale scenario}
%============================================================================

The second scenario to be evaluated is a large scale network of vehicles deployed around Maynooth town, based on the information
provided in Section \ref{S:Characteristics}. Here, we take into account a large number of cars (both moving and stationary), and also
a number of controlled intersections.

\begin{figure}[h]
\begin{centering}
\includegraphics[width=2.5in]{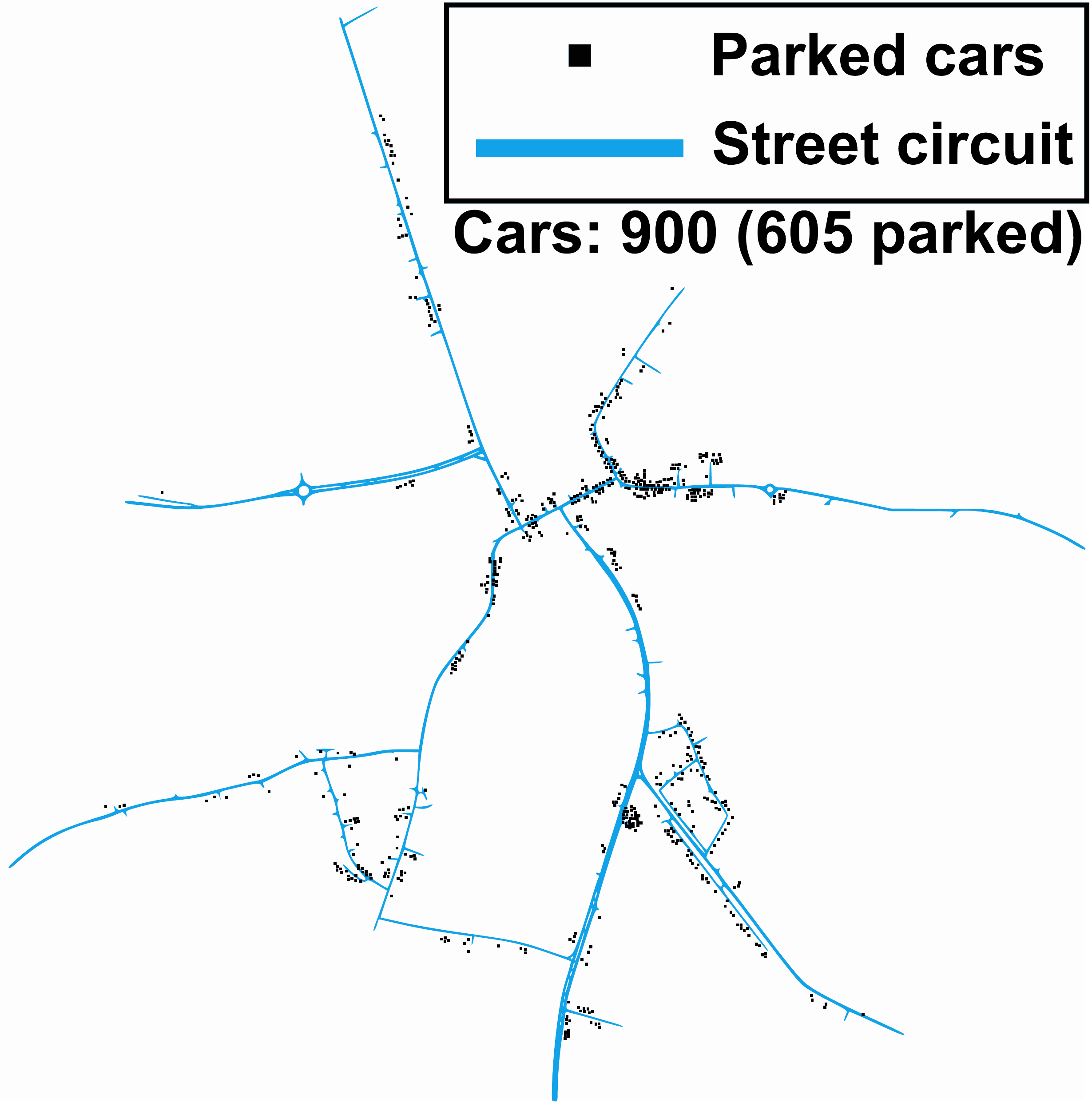}
\par\end{centering}
\caption{Chosen street circuit for the large scale test (blue) and selected parked cars (black).}
\label{fig:Street_circuit_Maynooth}
\end{figure}

The street circuit for this test is mostly composed by secondary and tertiary roads, and the parked cars to be used are those in a proximity of 15 m to the selected roads (a total of 605 parked cars), as shown in Fig. \ref{fig:Street_circuit_Maynooth}.
This assumption of 15 m for the communication range is consistent with a clear line-of-sight condition between parked cars and moving vehicles.

\begin{figure}[h]
\centering
\subfloat[Car 5, GCPSO with $\sigma_R = 0.2$ m.]{\includegraphics[width=3.5in]{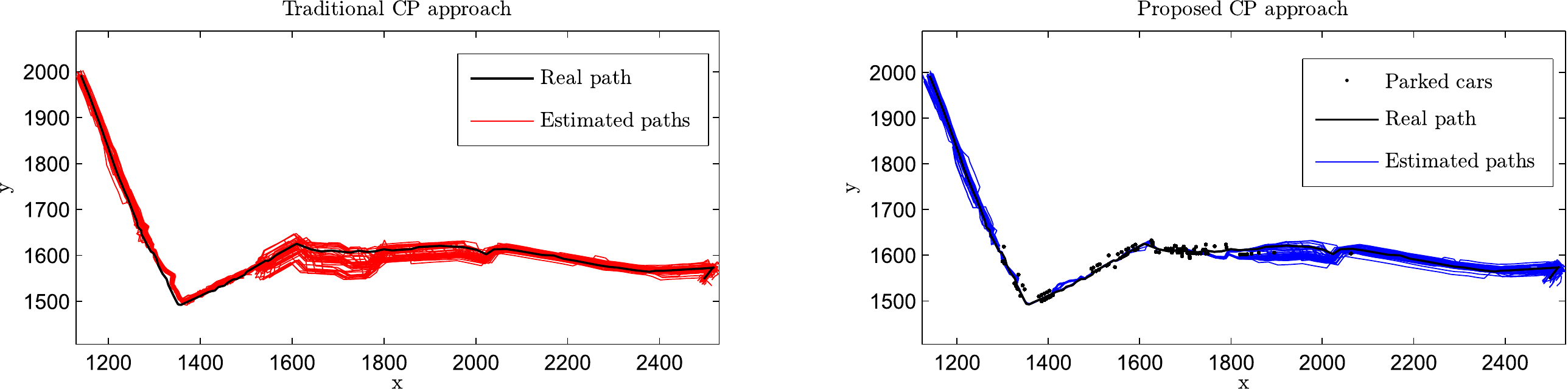}
\label{fig:Fig8a}}
\hfil
\\
\subfloat[Car 5, GCPSO with $\sigma_R = 4$ m.]{\includegraphics[width=3.5in]{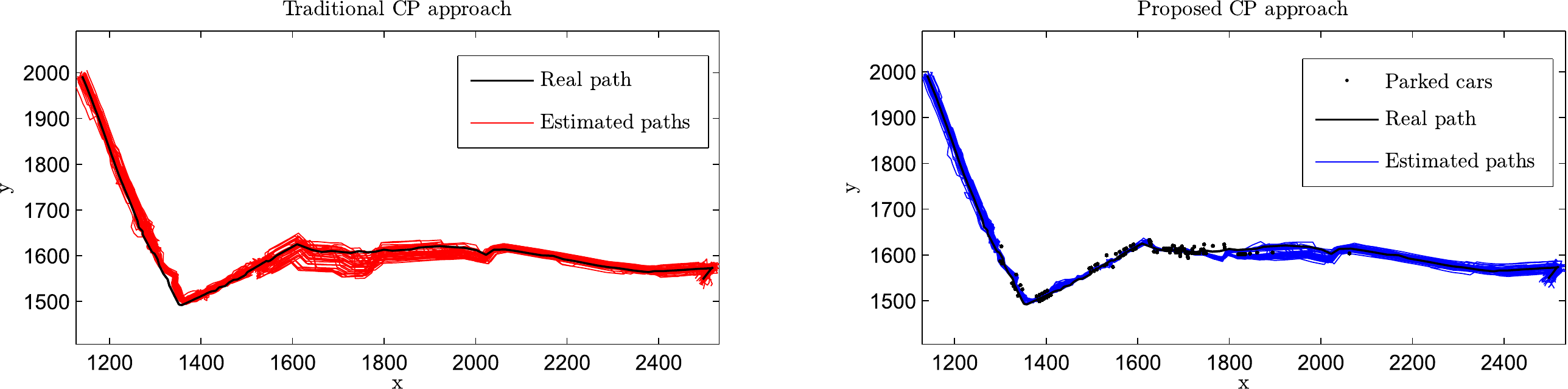}
\label{fig:Fig8b}}
\hfil
\\
\subfloat[Car 5, EKF with $\sigma_R = 0.2$ m.]{\includegraphics[width=3.5in]{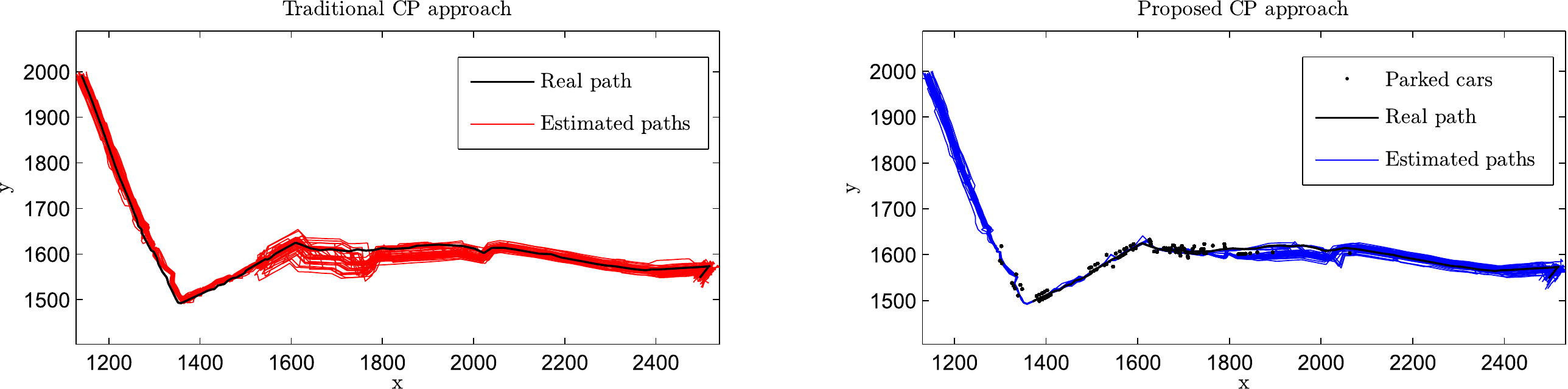}
\label{fig:Fig8c}}
\hfil
\\
\subfloat[Car 5, EKF with $\sigma_R = 4$ m.]{\includegraphics[width=3.5in]{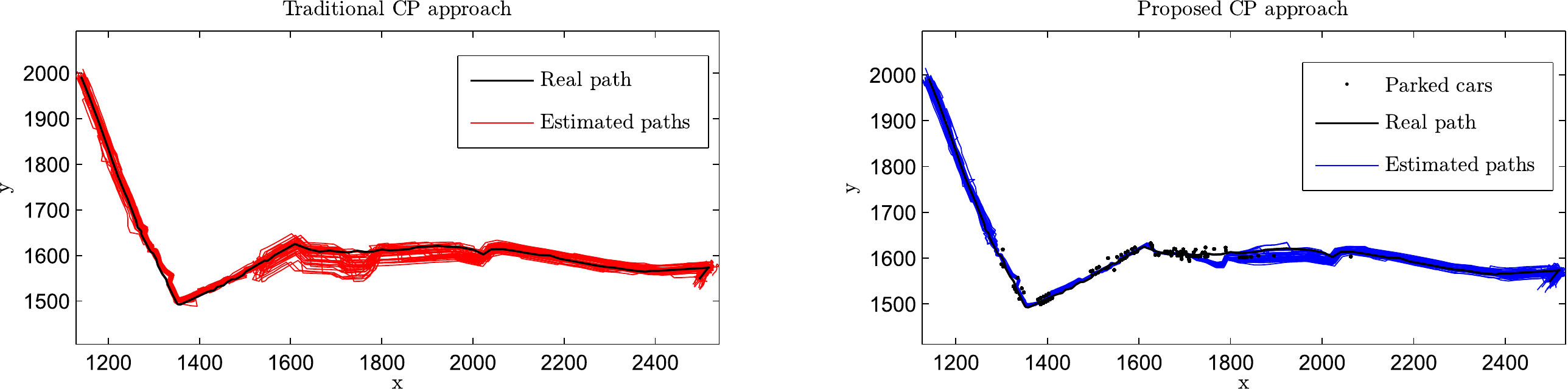}
\label{fig:Fig8d}}
\caption{Simulations results for Car 5 (40 different estimations).}
\label{fig:Results_Maynooth}
\end{figure}

Our experiments are constructed as follows. We use 295 blind moving cars are divided into two groups: 200 cars are selected with a random route and with initial position randomly distributed along the entire street circuit, and 95 other cars enter the street circuit at a rate of one every 4 seconds with a route defined by a random starting/ending pair of edges. The obtained results for the five cars with the longest routes are presented in Tables \ref{tab:RMSE_Example2_1}, \ref{tab:RMSE_Example2_2}, Fig. \ref{fig:Results_Maynooth} and Table \ref{tab:Parked_cars_per_Km}.

%=================================================================================================
% Table GCPSO
%=================================================================================================

\begin{table}
\caption{Results for the large scale scenario using GCPSO and 15 m of communication
zone: Localisation RMSE and improvements for the five cars with the
longest routes.}
\label{tab:RMSE_Example2_1}

\begin{centering}
\begin{tabular}{|c|c|c|c|c|c|c|}
\cline{3-6}
\multicolumn{1}{c}{} & & \multicolumn{4}{c|}{RMSE (meters)} & \multicolumn{1}{c}{}\tabularnewline
\cline{3-6}
\multicolumn{1}{c}{} & & \multicolumn{2}{c|}{Traditional} & \multicolumn{2}{c|}{Proposed} & \multicolumn{1}{c}{}\tabularnewline
\multicolumn{1}{c}{} & & \multicolumn{2}{c|}{CP approach} & \multicolumn{2}{c|}{CP approach} & \multicolumn{1}{c}{}\tabularnewline
\cline{3-7}
\multicolumn{1}{c}{} & & \multirow{2}{*}{Mean} & \multirow{2}{*}{$\sigma$} & \multirow{2}{*}{Mean} & \multirow{2}{*}{$\sigma$} & Average\tabularnewline
\multicolumn{1}{c}{} & & & & & & improvement\tabularnewline
\hline
\hline
\multirow{5}{*}{\begin{sideways}
$\sigma_{R}$= 0.2 m
\end{sideways}} & Car 1 & 8.30 & 7.14 & 6.92 & 7.11 & 16.59\%\tabularnewline
\cline{2-7}
& Car 2 & 8.63 & 5.23 & 6.71 & 5.71 & 22.20\%\tabularnewline
\cline{2-7}
& Car 3 & 9.44 & 7.26 & 5.79 & 6.87 & 38.71\%\tabularnewline
\cline{2-7}
& Car 4 & 9.80 & 7.16 & 5.48 & 5.57 & 44.05\%\tabularnewline
\cline{2-7}
& Car 5 & 10.67 & 7.20 & 3.96 & 5.22 & 62.86\%\tabularnewline
\hline
\hline
\multirow{5}{*}{\begin{sideways}
$\sigma_{R}$= 4 m
\end{sideways}} & Car 1 & 9.30 & 7.03 & 8.35 & 6.53 & 10.22\%\tabularnewline
\cline{2-7}
& Car 2 & 8.92 & 5.23 & 7.91 & 5.82 & 11.30\%\tabularnewline
\cline{2-7}
& Car 3 & 9.70 & 7.30 & 7.47 & 6.93 & 23.06\%\tabularnewline
\cline{2-7}
& Car 4 & 9.83 & 7.10 & 6.92 & 5.02 & 29.57\%\tabularnewline
\cline{2-7}
& Car 5 & 11.4 & 7.09 & 5.78 & 4.56 & 49.49\%\tabularnewline
\hline
\end{tabular}
\par\end{centering}

~

Mean: average of 40 different measurements.

$\sigma$: standard deviation.

\end{table}
%=================================================================================================

%=================================================================================================
% Table EKF
%=================================================================================================

\begin{table}
\caption{Results for the large scale scenario using EKF and 15 m of communication
zone: Localisation RMSE and improvements for the five cars with the
longest routes.}
\label{tab:RMSE_Example2_2}

\begin{centering}
\begin{tabular}{|c|c|c|c|c|c|c|}
\cline{3-6}
\multicolumn{1}{c}{} & & \multicolumn{4}{c|}{RMSE (meters)} & \multicolumn{1}{c}{}\tabularnewline
\cline{3-6}
\multicolumn{1}{c}{} & & \multicolumn{2}{c|}{Traditional} & \multicolumn{2}{c|}{Proposed} & \multicolumn{1}{c}{}\tabularnewline
\multicolumn{1}{c}{} & & \multicolumn{2}{c|}{CP approach} & \multicolumn{2}{c|}{CP approach} & \multicolumn{1}{c}{}\tabularnewline
\cline{3-7}
\multicolumn{1}{c}{} & & \multirow{2}{*}{Mean} & \multirow{2}{*}{$\sigma$} & \multirow{2}{*}{Mean} & \multirow{2}{*}{$\sigma$} & Average\tabularnewline
\multicolumn{1}{c}{} & & & & & & improvement\tabularnewline
\hline
\hline
\multirow{5}{*}{\begin{sideways}
$\sigma_{R}$= 0.2 m
\end{sideways}} & Car 1 & 7.77 & 6.67 & 7.22 & 7.25 & 7.06\%\tabularnewline
\cline{2-7}
& Car 2 & 7.44 & 4.20 & 6.69 & 5.62 & 10.19\%\tabularnewline
\cline{2-7}
& Car 3 & 9.12 & 7.75 & 5.99 & 7.11 & 34.32\%\tabularnewline
\cline{2-7}
& Car 4 & 9.17 & 6.88 & 5.74 & 6.17 & 37.37\%\tabularnewline
\cline{2-7}
& Car 5 & 10.84 & 6.76 & 3.65 & 5.12 & 66.35\%\tabularnewline
\hline
\hline
\multirow{5}{*}{\begin{sideways}
$\sigma_{R}$= 4 m
\end{sideways}} & Car 1 & 8.18 & 7.69 & 7.87 & 6.13 & 3.81\%\tabularnewline
\cline{2-7}
& Car 2 & 8.76 & 6.43 & 8.13 & 7.34 & 7.12\%\tabularnewline
\cline{2-7}
& Car 3 & 9.80 & 7.59 & 7.03 & 7.06 & 28.21\%\tabularnewline
\cline{2-7}
& Car 4 & 9.76 & 6.93 & 7.10 & 5.42 & 27.23\%\tabularnewline
\cline{2-7}
& Car 5 & 10.19 & 7.01 & 5.37 & 5.13 & 47.23\%\tabularnewline
\hline
\end{tabular}
\par\end{centering}

~

Mean: average of 40 different measurements.

$\sigma$: standard deviation.

\end{table}

%=================================================================================================

From Tables \ref{tab:RMSE_Example2_1} and \ref{tab:RMSE_Example2_2} it can be observed that the proposed approach to enhance CP gives an average improvement of 56.48\% in the localisation RMSE for the best case (with maximum at 66.35\% for EKF with $\sigma_R=0.2$ m, and minimum at 47.23\% for EKF with $\sigma_R=4$ m), with respect to the traditional approach.
In Table \ref{tab:Parked_cars_per_Km} it can be observed that the density of parked vehicles experienced by a given moving car affects considerably
the accuracy of the positioning. Note that despite the fact that radio ranging techniques cannot provide accurate inter-vehicle distance measurements for CP purposes,
simulation results have shown that their utility for the localisation of moving vehicles can be potentially greatly improved using the suggested
approach, even in the case of simplistic and imprecise CP algorithms.

In our experiments we used a short-range communication distance of 15 m, mainly for mitigating the
effects of shadowing and multipath (see references \cite{GiovanniBellusci2009,Gozalvez2012}).
Clearly, the results for the large scale scenario can be improved if an extended communication range is adopted
(as preliminary concluded in Section \ref{s:small_scale_example}), but, as we have already mentioned, in such a case the effects of shadowing and
multipath cannot be neglected.

\begin{table}[h]
\caption{Overall improvements of the localisation RMSE for the five cars with the
longest routes in function of the
number of parked cars experienced (*), travelled distance (**), and the ratio (*)/(**).}
\label{tab:Parked_cars_per_Km}

\begin{centering}
\begin{tabular}{|c|c|c|c|c|}
\cline{2-5} 
\multicolumn{1}{c|}{} & Parked cars & Travelled & Ratio & Overall \tabularnewline
\multicolumn{1}{c|}{} & experienced & path {[}Km{]} & ({*}) / ({*}{*}) & average\tabularnewline
\multicolumn{1}{c|}{} & ({*}) &  ({*}{*}) &  & improvement\tabularnewline
\hline 
Car 1 & 6 & 2.34 & 2.57 & 9.42\%\tabularnewline
\hline 
Car 2 & 18 & 2.62 & 6.86 & 12.70\%\tabularnewline
\hline 
Car 3 & 51 & 4.21 & 12.11 & 31.08\%\tabularnewline
\hline 
Car 4 & 58 & 5.10 & 11.37 & 34.56\%\tabularnewline
\hline 
Car 5 & 97 & 3.23 & 30.05 & 56.48\%\tabularnewline
\hline 
\end{tabular}
\par\end{centering}

\end{table}

%¶¶¶¶¶¶¶¶¶¶¶¶¶¶¶¶¶¶¶¶¶¶¶¶¶¶¶¶¶¶¶¶¶¶¶¶¶¶¶¶¶¶¶¶¶¶¶¶¶¶¶¶¶¶¶¶¶¶¶¶¶¶¶¶¶¶¶¶¶¶¶¶¶¶¶¶¶¶¶¶¶¶¶¶¶¶¶¶¶¶¶¶¶¶¶¶¶¶¶¶¶¶¶¶¶¶¶¶¶¶¶¶¶¶¶¶¶¶¶¶¶¶¶¶¶¶¶¶¶¶¶¶¶¶¶¶¶¶¶¶¶¶¶¶¶¶¶¶¶¶¶¶
%¶¶¶¶¶¶¶¶¶¶¶¶¶¶¶¶¶¶¶¶¶¶¶¶¶¶¶¶¶¶¶¶¶¶¶¶¶¶¶¶¶¶¶¶¶¶¶¶¶¶¶¶¶¶¶¶¶¶¶¶¶¶¶¶¶¶¶¶¶¶¶¶¶¶¶¶¶¶¶¶¶¶¶¶¶¶¶¶¶¶¶¶¶¶¶¶¶¶¶¶¶¶¶¶¶¶¶¶¶¶¶¶¶¶¶¶¶¶¶¶¶¶¶¶¶¶¶¶¶¶¶¶¶¶¶¶¶¶¶¶¶¶¶¶¶¶¶¶¶¶¶¶
\section{Conclusions and future work} \label{S:Conclusions}
%¶¶¶¶¶¶¶¶¶¶¶¶¶¶¶¶¶¶¶¶¶¶¶¶¶¶¶¶¶¶¶¶¶¶¶¶¶¶¶¶¶¶¶¶¶¶¶¶¶¶¶¶¶¶¶¶¶¶¶¶¶¶¶¶¶¶¶¶¶¶¶¶¶¶¶¶¶¶¶¶¶¶¶¶¶¶¶¶¶¶¶¶¶¶¶¶¶¶¶¶¶¶¶¶¶¶¶¶¶¶¶¶¶¶¶¶¶¶¶¶¶¶¶¶¶¶¶¶¶¶¶¶¶¶¶¶¶¶¶¶¶¶¶¶¶¶¶¶¶¶¶¶
%¶¶¶¶¶¶¶¶¶¶¶¶¶¶¶¶¶¶¶¶¶¶¶¶¶¶¶¶¶¶¶¶¶¶¶¶¶¶¶¶¶¶¶¶¶¶¶¶¶¶¶¶¶¶¶¶¶¶¶¶¶¶¶¶¶¶¶¶¶¶¶¶¶¶¶¶¶¶¶¶¶¶¶¶¶¶¶¶¶¶¶¶¶¶¶¶¶¶¶¶¶¶¶¶¶¶¶¶¶¶¶¶¶¶¶¶¶¶¶¶¶¶¶¶¶¶¶¶¶¶¶¶¶¶¶¶¶¶¶¶¶¶¶¶¶¶¶¶¶¶¶¶

CP is a localisation approach that is based on cooperation among several agents (nodes). In the case of VANETs, such agents can either be vehicles or road-side units (RSUs). Vehicles are traditionally the mobile part of the VANETs, while the RSUs constitute the static part. RSUs are, in general, fixed roadside stations that act like beacons or anchors depending on the accuracy with which their position is known (obtained via a dedicated calibration process). Despite the benefits of having a large number of RSUs in the VANETs, their deployment usually requires a complex and expensive process.

The contribution of this paper was, among others, to evaluate the use of on-street parked cars as an alternative to RSUs.
Such an alternative solution is generally cheaper, easier (no extra hardware is required), and more flexible than fixed infrastructure.
 Further, as the vehicle fleet refreshes itself after a number of years, technology updates
are easier to realise using vehicles as a platform. In addition, such new anchor nodes are usually
close to on-street blind cars, which is convenient to mitigate measurement noise due to shadowing and multipath effects of V2X measurements. Also,
the proximity of on-street stationary cars reduces the requirements on communication ranges, which in turn shortens the time required for accurate
localisation of blind cars. Finally, parked cars are dense, thus offering a potential platform for autonomous vehicles in urban situations.
Our simulations demonstrate the power of this technique in several settings. On-going work involves testing
these and related ideas  in cooperation with a vehicle manufacturer.

%¶¶¶¶¶¶¶¶¶¶¶¶¶¶¶¶¶¶¶¶¶¶¶¶¶¶¶¶¶¶¶¶¶¶¶¶¶¶¶¶¶¶¶¶¶¶¶¶¶¶¶¶¶¶¶¶¶¶¶¶¶¶¶¶¶¶¶¶¶¶¶¶¶¶¶¶¶¶¶¶¶¶¶¶¶¶¶¶¶¶¶¶¶¶¶¶¶¶¶¶¶¶¶¶¶¶¶¶¶¶¶¶¶¶¶¶¶¶¶¶¶¶¶¶¶¶¶¶¶¶¶¶¶¶¶¶¶¶¶¶¶¶¶¶¶¶¶¶¶¶¶¶
%¶¶¶¶¶¶¶¶¶¶¶¶¶¶¶¶¶¶¶¶¶¶¶¶¶¶¶¶¶¶¶¶¶¶¶¶¶¶¶¶¶¶¶¶¶¶¶¶¶¶¶¶¶¶¶¶¶¶¶¶¶¶¶¶¶¶¶¶¶¶¶¶¶¶¶¶¶¶¶¶¶¶¶¶¶¶¶¶¶¶¶¶¶¶¶¶¶¶¶¶¶¶¶¶¶¶¶¶¶¶¶¶¶¶¶¶¶¶¶¶¶¶¶¶¶¶¶¶¶¶¶¶¶¶¶¶¶¶¶¶¶¶¶¶¶¶¶¶¶¶¶¶
\section*{Acknowledgment}
%¶¶¶¶¶¶¶¶¶¶¶¶¶¶¶¶¶¶¶¶¶¶¶¶¶¶¶¶¶¶¶¶¶¶¶¶¶¶¶¶¶¶¶¶¶¶¶¶¶¶¶¶¶¶¶¶¶¶¶¶¶¶¶¶¶¶¶¶¶¶¶¶¶¶¶¶¶¶¶¶¶¶¶¶¶¶¶¶¶¶¶¶¶¶¶¶¶¶¶¶¶¶¶¶¶¶¶¶¶¶¶¶¶¶¶¶¶¶¶¶¶¶¶¶¶¶¶¶¶¶¶¶¶¶¶¶¶¶¶¶¶¶¶¶¶¶¶¶¶¶¶¶
%¶¶¶¶¶¶¶¶¶¶¶¶¶¶¶¶¶¶¶¶¶¶¶¶¶¶¶¶¶¶¶¶¶¶¶¶¶¶¶¶¶¶¶¶¶¶¶¶¶¶¶¶¶¶¶¶¶¶¶¶¶¶¶¶¶¶¶¶¶¶¶¶¶¶¶¶¶¶¶¶¶¶¶¶¶¶¶¶¶¶¶¶¶¶¶¶¶¶¶¶¶¶¶¶¶¶¶¶¶¶¶¶¶¶¶¶¶¶¶¶¶¶¶¶¶¶¶¶¶¶¶¶¶¶¶¶¶¶¶¶¶¶¶¶¶¶¶¶¶¶¶¶

This work was supported in part by both the Science Foundation Ireland grant 11/PI/1177, and the European Commission under TEAM, a large scale integrated project part of the FP7-ICT for Cooperative Systems for energy-efficient and sustainable mobility.

\bibliographystyle{ieeetr}
\bibliography{Biblio_Parked}

\begin{IEEEbiography}[{\includegraphics[width=1in,height=1.25in,clip,keepaspectratio]{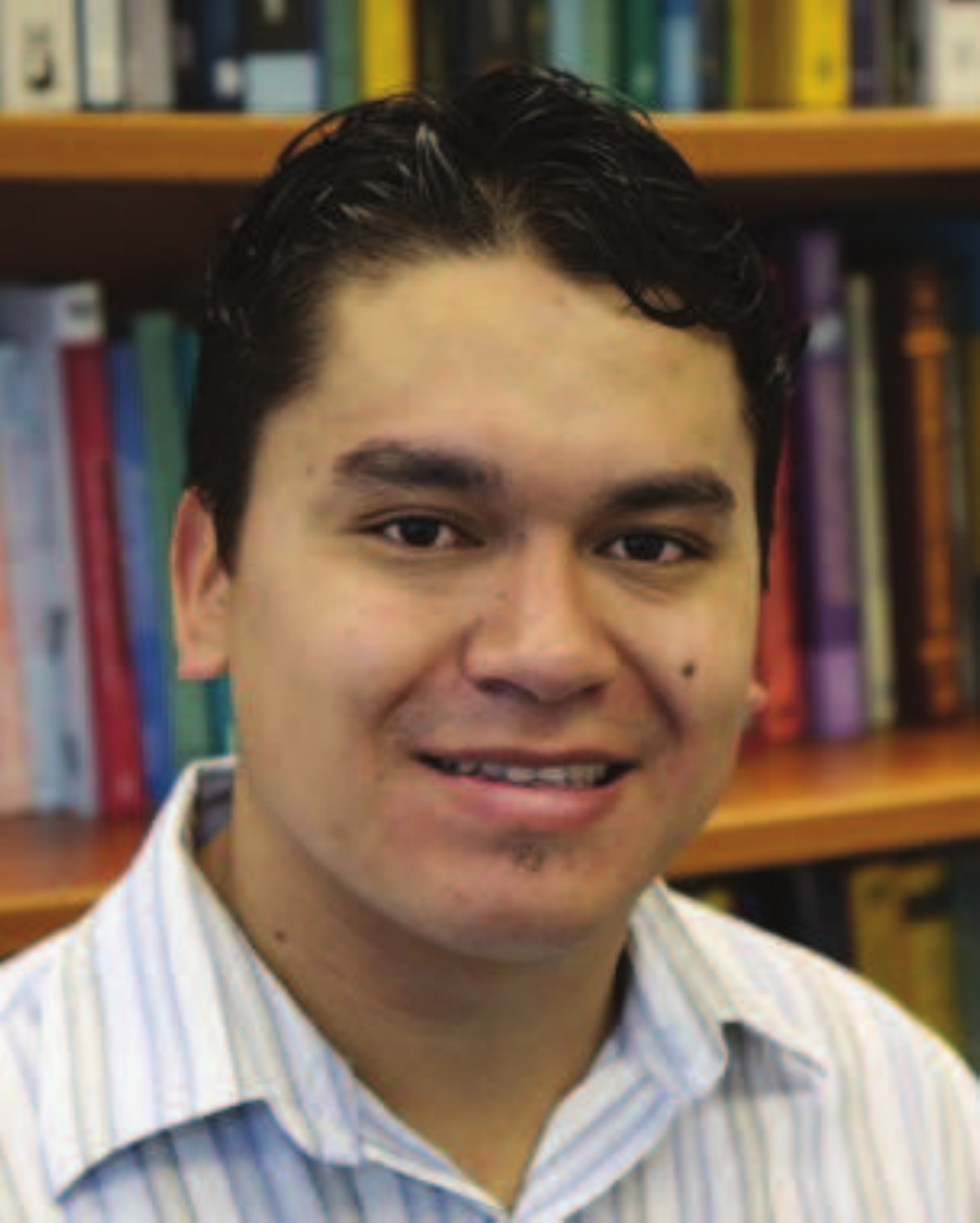}}]{Rodrigo H. Ord\'{o}\~{n}ez-Hurtado}
received his degree in Engineering in Industrial Automatica from the University of Cauca, Colombia, in 2005. He commenced his Ph.D. in Electrical Engineering at the University of Chile in March 2008, and between April-November 2012, he was an intern at the Hamilton Institute, a multi-disciplinary research centre established at the National University of Ireland in Maynooth. In November 2012, he received his Ph.D. degree from the University of Chile, and since December 2012 he has held a postdoctoral position at the Hamilton Institute, working with Professor R. Shorten and his research group. Rodrigo's interests include robust adaptive systems (control and identification), stability of switched systems, swarm intelligence, large-scale systems and intelligent transportation systems. His focus is on applications to the mining industry and transportation systems.
\end{IEEEbiography}

\begin{IEEEbiography}[{\includegraphics[width=1in,height=1.25in,clip,keepaspectratio]{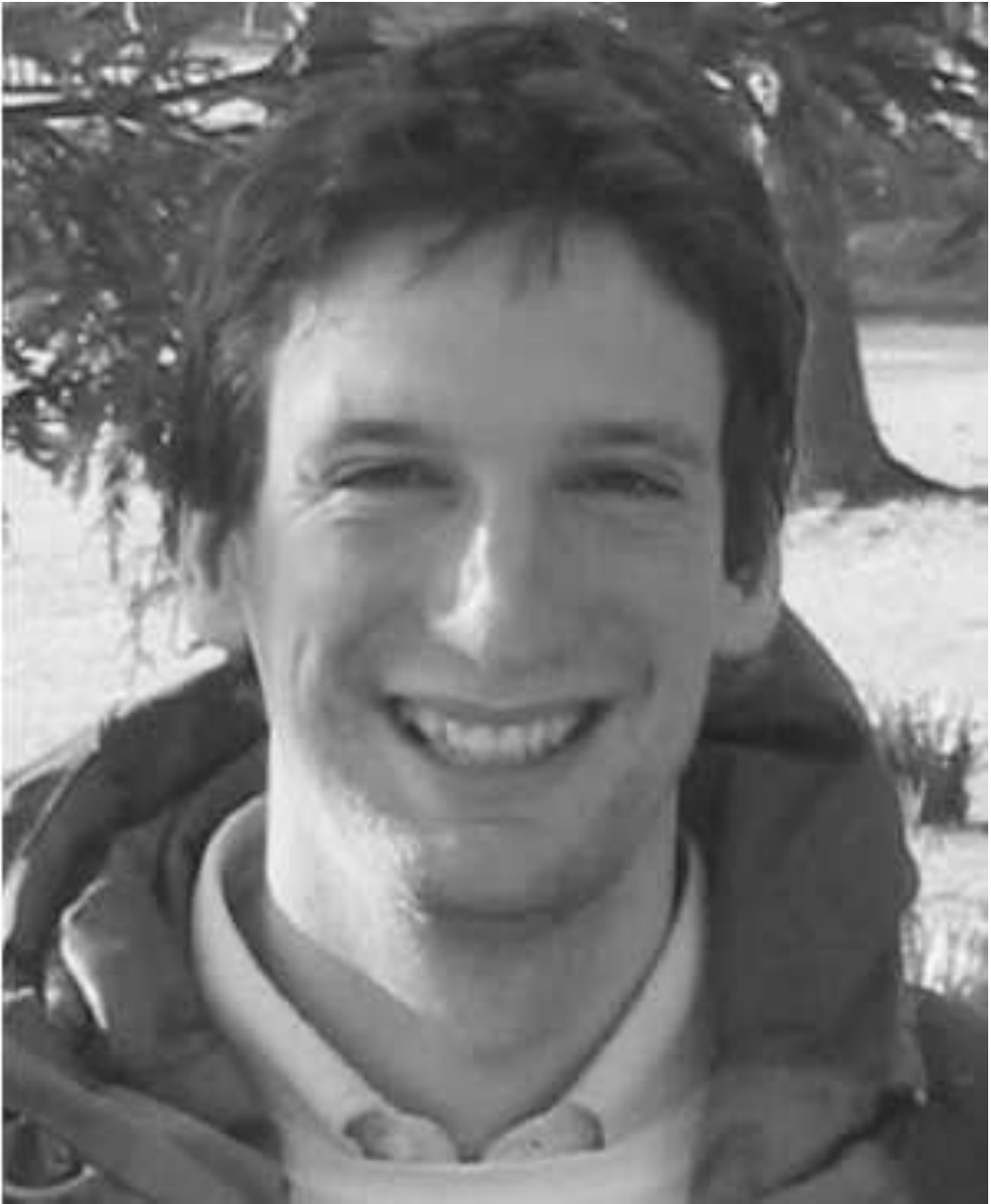}}]{Emanuele Crisostomi}
received the B.S. degree in computer science engineering, the M.S. degree in automatic control, and the Ph.D. degree in automatics, robotics, and bioengineering, from the University of Pisa, Italy, in 2002, 2005, and 2009, respectively. He is currently an Assistant Professor of electrical engineering with the Department of Energy, Systems, Territory and Constructions Engineering, University of Pisa. His research interests include control and optimization of large-scale systems, with applications to smart grids and green mobility networks.
\end{IEEEbiography}

\begin{IEEEbiography}[{\includegraphics[width=1in,height=1.25in,clip,keepaspectratio]{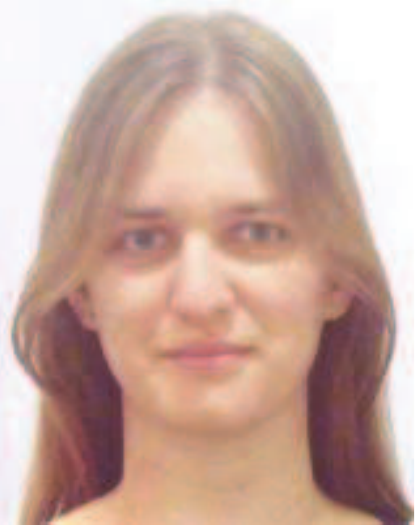}}]{Wynita Griggs}
received her B.Sc. (Hons) degree in Mathematics from the University of Queensland, Brisbane, Australia in 2002 and her Ph.D. in Engineering from the Australian National University, Canberra, Australia in 2007. She is currently a postdoctoral research fellow at the Hamilton Institute, National University of Ireland, Maynooth. Her research interests include stability theory with applications to feedback control systems; and smart transportation.
\end{IEEEbiography}

\begin{IEEEbiography}[{\includegraphics[width=1in,height=1.25in,clip,keepaspectratio]{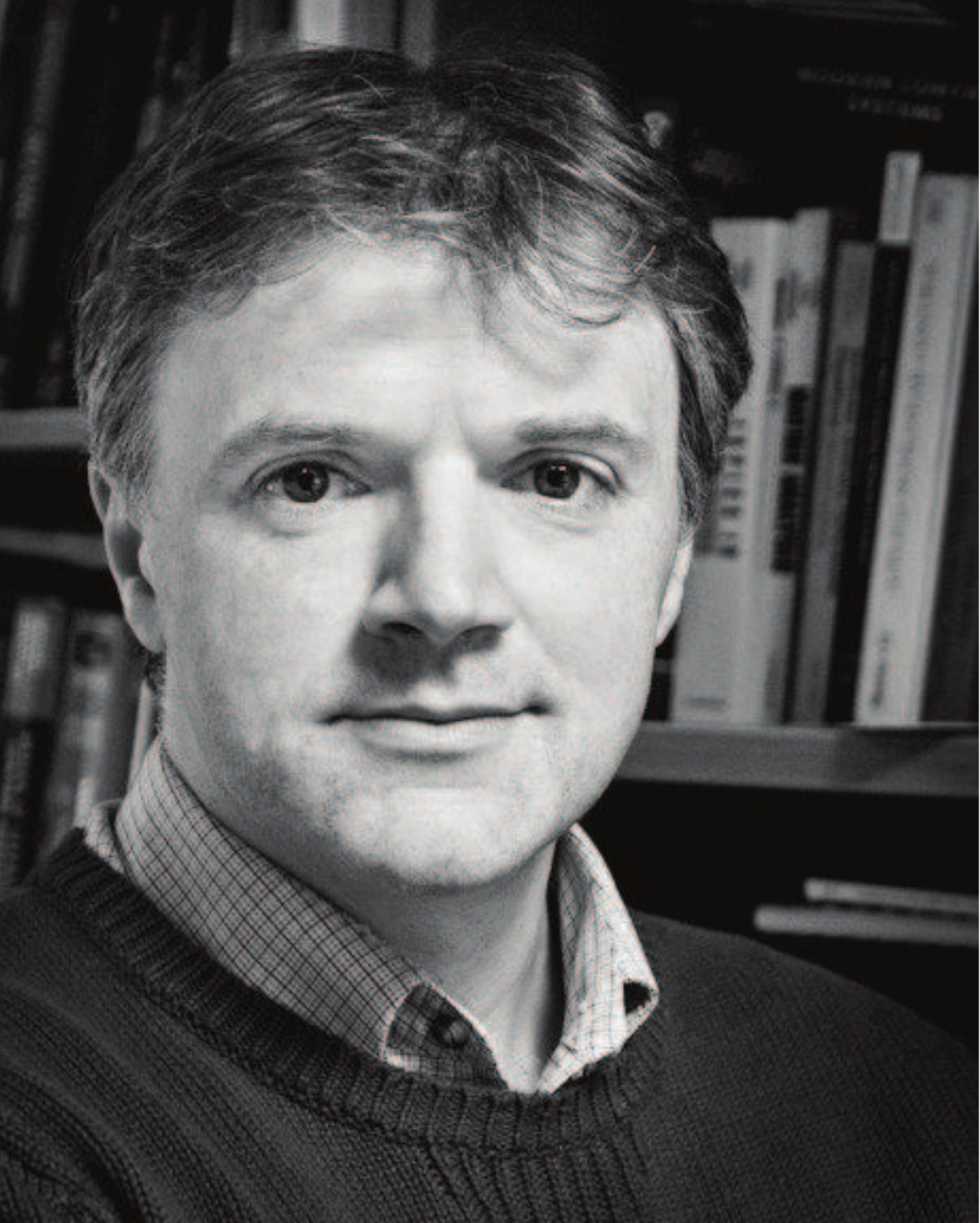}}]{Robert Shorten}
graduated from University College Dublin (UCD) in 1990 with a First Class Honours B.E. degree in Electronic Engineering. He was awarded a Ph.D. in 1996, also from UCD, while based at Daimler-Benz Research in Berlin, Germany. From 1993 to 1996, Prof. Shorten was the holder of a Marie Curie Fellowship at Daimler-Benz Research to conduct research in the area of smart gearbox systems. Following a brief spell at the Center for Systems Science, Yale University, working with Professor K. S. Narendra, Prof. Shorten returned to Ireland as the holder of a European Presidency Fellowship in 1997. Prof. Shorten is a co-founder of the Hamilton Institute at NUI Maynooth, where he was a full Professor until March 2013. He was also a Visiting Professor at TU Berlin from 2011-2012. Professor Shorten is currently a senior research manager at IBM Research Ireland. Prof. Shorten's research spans a number of areas. He has been active in computer networking, automotive research, collaborative mobility (including smart transportation and electric vehicles), as well as basic control theory and linear algebra. His main field of theoretical research has been the study of hybrid dynamical systems.
\end{IEEEbiography}

\end{document}